\documentclass[11pt]{article}

\usepackage{amssymb, amsfonts, amsmath, amsthm, amscd, graphicx, color, textcomp}
\usepackage[colorlinks=true,breaklinks=true,bookmarks=true,urlcolor=blue,
citecolor=blue,linkcolor=blue,bookmarksopen=false,draft=false]{hyperref}
\usepackage[symbol]{footmisc}
\usepackage{mathrsfs}
\usepackage{cite}
\usepackage{secdot}
\usepackage{relsize}	
\newtheorem{theorem}{Theorem}[section]
\newtheorem{corollary}[theorem]{Corollary}
\newtheorem{lemma}[theorem]{Lemma}
\newtheorem{proposition}[theorem]{Proposition}
\theoremstyle{definition}
\newtheorem{definition}[theorem]{Definition}
\newtheorem{example}[theorem]{Example}

\begin{document}

\begin{center}
{\Large \textbf {Existence and uniqueness of fixed point on closed ball in  multiplicative $ \mathbf{G}_\mathcal{M}-$metric space}}\medskip \smallskip

\textbf{Mohamed Gamal}$^{1,}$\footnote{Corresponding author: m$\_$gamal29@sci.svu.edu.eg (M. Gamal),},
\textbf{Fu-Gui Shi}$^{2}$\medskip

\noindent$^{1}$Department of Mathematics, Faculty of Science, South Valley University, Qena 83523, Egypt\\
\noindent$^{2}$School of Mathematics and Statistics, Beijing Institute of Technology, Beijing 102488, China\\
\medskip

\noindent Email: m$\_$gamal29@sci.svu.edu.eg,$\:$ fuguishi@bit.edu.cn

\noindent\hrulefill
\end{center}

%%%%%%%%%%%%%%%%%%%%%%%%%%%%%%%%%%%%%%%%%%%%%%%%%%%%%%%%%%%%%%%%%%%%%%%%%%%%%%%%%%%%%%%%%%
\noindent{\large{\textbf{Abstract}}\smallskip

In this article, we investigate some fixed point results satisfying a new generalized $\Delta$-implicit contractive condition in ordered complete multiplicative $ \mathbf{G}_\mathcal{M}-$metric space.  Also, some new definitions and fixed point theorems are presented in ordered complete multiplicative $ \mathbf{G}_\mathcal{M}-$metric space. Furthermore, some nontrivial and illustrative examples are given to validate our obtained results.\smallskip

%%%%%%%%%%%%%%%%%%%%%%%%%%%%%%%%%%%%%%%%%%%%%%%%%%%%%%%%%%%%%%%%%%%%%%%%%%%%%%%%EE%%%%%%%%%%
\noindent{\large{\textbf{Keywords:}}
Multiplicative $ \mathbf{G}_\mathcal{M}-$metric, closed ball, generalized $\Delta$-implicit contraction, fixed point methodology. \smallskip 

%%%%%%%%%%%%%%%%%%%%%%%%%%%%%%%%%%%%%%%%%%%%%%%%%%%%%%%%%%%%%%%%%%%%%%%%%%%%%%%%%%%%%%%%%%%%
\smallskip

\noindent\hrulefill

%%%%%%%%%%%%%%%%%%%%%%%%%%%%%%%%%%%%%%%%%%%%%%%%%%%%%%%%%%%%%%%%%%%%%%%%%%%%%%%%%%%%%%%%%%%%
\section{\bf{Introduction}}

\indent \indent Fixed point theory ($\mathbf{FPT}$) plays a vivid, exciting and fundamental role in the field of functional analysis. S. Banach \cite{Banach1992} presented a foundational principle and it becomes a vital tool in the field of metric fixed point with a lot of applications to ensure the existence and uniqueness of fixed point ($\mathbf{FP}$). This theorem is called  Banach fixed-point theorem (also is known as contractive mapping theorem or  contraction mapping theorem). Due to its advantages, many authors showed different improvements and extensions of this theorem in various distance spaces (see \cite{Abbas, Ali, Banach1992, Bashirov1, Bhatt, Bojor, Boriceanu, Butnariu, Debnath, Dosenovic1, Dosenovic, Gu1,  Hussain, Jachymski, Jain, Piao, He, Jiang, Kang1, Todorcevic}). \medskip

Bashirov et al. \cite{Bashirov1} established the concept of multiplicative calculus, showed the foundational theorem of multiplicative calculus and studied some fundamental properties.  After that, other properties in multiplicative metric space ($\:M^{\circ}M^{\bullet}S$) were studied and constructed in \cite{Bashirov2,Florack}. In 2012, \"{O}zavsar et al. \rmfamily\cite{Ozavsar} displayed the definition of multiplicative contraction mappings on $\:M^{\circ}M^{\bullet}S \:$ in such a method that multiplicative triangle inequality is used instead of the usual triangular inequality and presented different existence results of $\: \mathbf{FP} \:$ beside various topological characteristics of $\:M^{\circ}M^{\bullet}S$. Also, many researchers studied fixed point theorems in multiplicative metric space using  weak commutative mappings, locally contractive mappings, $\mathcal{E.A}$-property, compatible-type mappings and  generalized contraction mappings with cyclic ($\alpha,\beta$)-admissible mapping respectively, for more illustrations (see \cite{He,Abbas, Abdou,Srinivas,Yamaod}). In 2016, Nagpal et al. \cite{Nagpal} introduced the concept of multiplicative generalized metric space and studied the notion of weakly commuting, compatible maps and its variants, weakly compatible, weakly compatible with properties $ (\mathcal{CLR}) $ and ($\mathcal{E.A}$) in the same space. \medskip 

According to this orientation, the major purpose of this article is to prvoe some new fixed point theorems satisfying a new generalized $\Delta$-implicit contraction on a closed ball in ordered complete multiplicative $ \mathbf{G}_\mathcal{M}-$metric space ($M^{\circ}\, \mathbf{G}_\mathcal{M}-M^{\bullet}S$). Eventually, we prove some nontrivial examples to support new results.

%%%%%%%%%%%%%%%%%%%%%%%%%%%%%%%%%%%%%%%%%%%%%%
\section{\bf{Preliminary}}

\indent \indent Now, we recall some well-known notations and definitions that will be used in our subsequent discussion.

%%%%%%%%%%%%%%%%%%%%%%%%%%%%%%%%%%%%%%%%%%%%%%
\begin{definition}
	
{\rm\cite{Banach1992}} Let $\: \mathcal{P} \:$ be a mapping on a nonempty set $\: \mathcal{Z} \:$. Then a point $\check{\nu} \in \mathcal{Z} $ is called a $\:\mathbf{FP} \:$ of $\: \mathcal{P} \:$ if $\: \mathcal{P}\check{\nu} = \check{\nu}. $
	
\end{definition}

%%%%%%%%%%%%%%%%%%%%%%%%%%%%%%%%%%%%%%%%%%%%%%
\begin{definition} 
	
{\rm\cite{Bashirov1}} %Bashirov 2008
Let a non-empty set $\: \mathcal{Z} \:$ and $\: \zeta_\mathcal{M}: \mathcal{Z} \times \mathcal{Z} \longrightarrow \mathbb{R^+}$ be a function satisfying the following properties:
	
\noindent ($ \zeta_1 $) $\: \zeta_\mathcal{M}(\check{\nu},\check{\varpi}) \: \geq \: 1, \quad$  $\forall \:\:\: \check{\nu},\check{\varpi} \in \mathcal{Z} $;
	
\noindent ($ \zeta_2 $) $\: \zeta_\mathcal{M}(\check{\nu},\check{\varpi}) \: =  \: 1 \quad \:$ iff $\:\:  \check{\nu}=\check{\varpi} $;
	
\noindent ($ \zeta_3 $) $\: \zeta_\mathcal{M}(\check{\nu},\check{\varpi}) \: = \: \zeta_\mathcal{M}(\check{\varpi},\check{\nu}) \qquad$ (symmetry);
	
\noindent ($\zeta_4$) $\: \zeta_\mathcal{M}(\check{\nu},\check{\varpi}) \:\: \leq \:  \zeta_\mathcal{M}(\check{\nu},\check{\vartheta}) \: . \: \zeta_\mathcal{M}(\check{\vartheta},\check{\varpi}) \qquad \forall \:\:\: \check{\nu},\check{\varpi},\check{\vartheta} \in \mathcal{Z} \quad$ (multiplicative triangle inequality).\smallskip
	
\noindent Then, $ \zeta_\mathcal{M} $ is  a multiplicative metric on $ \mathcal{Z} $ and the pair $ (\mathcal{Z},\zeta_\mathcal{M})$ is a $\:M^{\circ}M^{\bullet}S$.
	
\end{definition} 

%%%%%%%%%%%%%%%%%%%%%%%%%%%%%%%%%%%%%%%%%%%%%
\begin{definition} 
	
{\rm\cite{Nagpal}} 
Suppose $\: \mathcal{Z} \:$ be a non-empty set and $\: {\mathbf{G}_\mathcal{M}}: \mathcal{Z} ^3  \longrightarrow \mathbb{R^+} $ be a function satisfying the following conditions:
	
\noindent ($\mathbf{G}_{\mathcal{M}_1}$) $\: {\mathbf{G}_\mathcal{M}}(\check{\nu},\check{\varpi},\check{\vartheta}) \: = \: 1 \quad$ if $\: \check{\nu} = \check{\varpi} = \check{\vartheta}$;
	
\noindent ($\mathbf{G}_{\mathcal{M}_2}$) $\: 1 \: < \: {\mathbf{G}_\mathcal{M}}(\check{\nu},\check{\nu},\check{\varpi})  \qquad  \forall \:\:\: \check{\nu},\check{\varpi} \in  \mathcal{Z}  \:$ with $\: \check{\nu} \neq \check{\varpi} $;
	
\noindent ($\mathbf{G}_{\mathcal{M}_3}$)  $\: {\mathbf{G}_\mathcal{M}}(\check{\nu},\check{\nu},\check{\varpi}) \: \leq \:  {\mathbf{G}_\mathcal{M}}(\check{\nu},\check{\varpi},\check{\vartheta}) \qquad \forall \:\:\: \check{\nu},\check{\varpi},\check{\vartheta} \in \mathcal{Z}  \:\:$ with $\:\: \check{\varpi} \neq \check{\vartheta}$;
	
\noindent ($\mathbf{G}_{\mathcal{M}_4}$) $\: {\mathbf{G}_\mathcal{M}}(\check{\nu},\check{\varpi},\check{\vartheta}) \: = \: {\mathbf{G}_\mathcal{M}}(\check{\nu},\check{\vartheta},\check{\varpi}) \: = \: {\mathbf{G}_\mathcal{M}}(\check{\varpi},\check{\vartheta},\check{\nu}) \: = \: ... \qquad $ (symmetry);
	
\noindent ($\mathbf{G}_{\mathcal{M}_5}$) $\: {\mathbf{G}_\mathcal{M}}(\check{\nu},\check{\varpi},\check{\vartheta}) \leq {\mathbf{G}_\mathcal{M}}(\check{\nu},\check{\tau},\check{\tau}). {\mathbf{G}_\mathcal{M}}(\check{\tau},\check{\varpi},\check{\vartheta}) \quad \forall \:\: \check{\nu},\check{\varpi},\check{\vartheta}, \check{\tau} \in \mathcal{Z},\: $ (rectangular inequality).\smallskip

\noindent Then, the function $ {\mathbf{G}_\mathcal{M}} $ is called a multiplicative generalized metric or, more accurately,  multiplicative $ \mathbf{G}_\mathcal{M}-$metric on $ \mathcal{Z}  $ and the pair ($\mathcal{Z}, \mathbf{G}_\mathcal{M} $) is a  $\:M^{\circ}\, \mathbf{G}_\mathcal{M}-M^{\bullet}S$.
	
\end{definition}

%%%%%%%%%%%%%%%%%%%%%%%%%%%%%%%%%%%%%%%%%%%%%%
\begin{definition}
Suppose  $ (\mathcal{Z}, \mathbf{G}_\mathcal{M}) $ be a $\, M^{\circ}\, \mathbf{G}_\mathcal{M}-M^{\bullet}S\,$ then for $ \check{\nu}_{0} \in \mathcal{Z}$, the $\: \mathbf{G}_\mathcal{M}-$ball with centre $\: \check{\nu}_{0} \:$ and radius $\: \gamma \:\: (\gamma >0) $ is,
$$\overline{\mathcal{B}_{\gamma} (\check{\nu}_{0}, \gamma)} = \{ \varrho \in \mathcal{Z}: \mathbf{G}_\mathcal{M}\: (\check{\nu}_{0}, \varrho, \varrho ) \leq \gamma \}.$$
\end{definition}

%%%%%%%%%%%%%%%%%%%%%%%%%%%%%%%%%%%%%%%%%%%%%

\begin{example}
%example 2.3, 2.4 in paper and Exxample 2.5 in another paper
\noindent Let $ (\mathcal{Z} ,d) $ be a $\: M^{\circ} M^{\bullet}S \:$ and $ {\mathbf{G}_\mathcal{M}}:\mathcal{Z}^3 \longrightarrow \mathbb{R^+} $ is defined by ${\mathbf{G}_\mathcal{M}}(\check{\nu},\check{\varpi},\check{\vartheta})=d(\check{\nu},\check{\varpi})\:.\:d(\check{\varpi},\check{\vartheta})\:.\:d(\check{\vartheta},\check{\nu})\:\:\: \forall \:\:  \check{\nu},\check{\varpi},\check{\vartheta} \in \mathcal{Z}. $ Then, $ {\mathbf{G}_\mathcal{M}}  $ is a multiplicative $ \mathbf{G}_\mathcal{M}-$metric on $ \mathcal{Z}  $ and $ (\mathcal{Z} ,{\mathbf{G}_\mathcal{M}}) $ is called $\:M^{\circ}\, \mathbf{G}_\mathcal{M}-M^{\bullet}S$.

\end{example}

%%%%%%%%%%%%%%%%%%%%%%%%%%%%%%%%%%%%%%%%%%%%%%
\begin{example}
	
\noindent Assume that $ (\mathcal{Z},d) $ be a usual metric space and $ {\mathbf{G}_\mathcal{M}}:\mathcal{Z}^{3} \longrightarrow \mathbb{R}^+ $ is defined by ${\mathbf{G}_\mathcal{M}}(\check{\nu},\check{\varpi},\check{\vartheta})=e\:^{d(\check{\nu},\check{\varpi}) \: + \: d(\check{\varpi},\check{\vartheta}) \: + \: d(\check{\vartheta}, \check{\nu})}\quad \forall \:\:  \check{\nu},\check{\varpi},\check{\vartheta} \in \mathcal{Z}. $ Thus, $ {\mathbf{G}_\mathcal{M}}  $ is a multiplicative $ \mathbf{G}_\mathcal{M}-$metric on $ \mathcal{Z} $ and $ (\mathcal{Z},{\mathbf{G}_\mathcal{M}}) $ is called $\:M^{\circ}\, \mathbf{G}_\mathcal{M}-M^{\bullet}S$.
	
\end{example}

%%%%%%%%%%%%%%%%%%%%%%%%%%%%%%%%%%%%%%%%%%%%%%

\begin{definition}\cite{Abbas1} Suppose  $ (\mathcal{Z}, \leq ) $  be a poset. Then, $ \check{\varpi}, \check{\rho} \in \mathcal{Z} $ are called comparable if $\: \check{\varpi} \leq \check{\rho} \:$ or $\: \check{\rho}  \leq \check{\varpi} \:$ holds.

\end{definition}

%%%%%%%%%%%%%%%%%%%%%%%%%%%%%%%%%%%%%%%%%%%%%

\begin{proposition} 
{\rm\cite{Nagpal}} Let $(\mathcal{Z},{\mathbf{G}_\mathcal{M}})$ be a $\:M^{\circ}\, \mathbf{G}_\mathcal{M}-M^{\bullet}S$. Then, for all $\: \check{\nu},\check{\varpi},\check{\vartheta},\check{\tau} \in \mathcal{Z}, $ the following properties are satisfying:\smallskip
	
\noindent $ (1) $ $\: {\mathbf{G}_\mathcal{M}}(\check{\nu},\check{\varpi},\check{\vartheta}) \: = \: 1 \quad$ if $\:\: \check{\nu} = \check{\varpi} = \check{\vartheta} $;\\
%%%%%%%%%%%%%%%%%%%%%%%%%%%%%%%%%%%%%%%%%%%%%%
$ (2) $ $\: {\mathbf{G}_\mathcal{M}}(\check{\nu},\check{\varpi},\check{\vartheta}) \: \leq \: {\mathbf{G}_\mathcal{M}}(\check{\nu},\check{\tau},\check{\tau})\: . \: {\mathbf{G}_\mathcal{M}}(\check{\varpi},\check{\tau},\check{\tau})\: . \: {\mathbf{G}_\mathcal{M}}(\check{\vartheta},\check{\tau},\check{\tau})$;\\
%%%%%%%%%%%%%%%%%%%%%%%%%%%%%%%%%%%%%%%%%%%%%%
$ (3) $ $\: {\mathbf{G}_\mathcal{M}}(\check{\nu},\check{\varpi},\check{\vartheta}) \: \leq \: {\mathbf{G}_\mathcal{M}}(\check{\nu},\check{\nu},\check{\varpi})\: . \: {\mathbf{G}_\mathcal{M}}(\check{\nu},\check{\nu},\check{\vartheta});$\\
%%%%%%%%%%%%%%%%%%%%%%%%%%%%%%%%%%%%%%%%%%%%%%
$ (4) $ $\: {\mathbf{G}_\mathcal{M}}(\check{\nu},\check{\varpi},\check{\varpi}) \: \leq \:  \mathbf{G}^{2}_\mathcal{M}(\check{\varpi},\check{\nu},\check{\nu})  $. 
	
\end{proposition}

%%%%%%%%%%%%%%%%%%%%%%%%%%%%%%%%%%%%%%%%%%%%%
\begin{lemma} 
	
{\rm \cite{Nagpal}} Let $ \{\check{\nu}_{k}\} $ be a sequence in a $\: M^{\circ}\, \mathbf{G}_\mathcal{M}-M^{\bullet}S$ $(\mathcal{Z},{\mathbf{G}_\mathcal{M}})$. If the sequence $ \{\check{\nu}_{k}\} $ is multiplicative $\mathbf{G}_\mathcal{M}-$convergent then it is multiplicative $\mathbf{G}_\mathcal{M}-$Cauchy ($ M^{\circ}\,\mathbf{G}_\mathcal{M}-C^{\bullet}$) sequence. 
	
\end{lemma}

%%%%%%%%%%%%%%%%%%%%%%%%%%%%%%%%%%%%%%%%%%%%%
\begin{lemma} {\rm\cite{Nagpal}} Let $ \{\check{\nu}_{k}\} $ be a sequence in a $\:M^{\circ}\, \mathbf{G}_\mathcal{M}-M^{\bullet}S$ $(\mathcal{Z},{\mathbf{G}_\mathcal{M}})$. The sequence $ \{\check{\nu}_{k}\} $ in $ \mathcal{Z} $ is multiplicative $\mathbf{G}_\mathcal{M}-$convergent to $ p \in \mathcal{Z}$ iff $\: {\mathbf{G}_\mathcal{M}}(\check{\nu}_{k},p,p)\longrightarrow 1, \:\:$ as $\:\: k \longrightarrow +\infty$.  
	
\end{lemma}

%%%%%%%%%%%%%%%%%%%%%%%%%%%%%%%%%%%%%%%%%%%%%%%%%%%%%%%%%%%%%%%%%%%%%%%%%%%%%%%%%%%%%%%%%%%%
\section{Main Results}

\indent \indent Now, we present our main theorem in ordered complete multiplicative $ \mathbf{G}_\mathcal{M}-$metric space.
%%%%%%%%%%%%%%%%%%%%%%%%%%%%%%%%%%%%%%%%%%%%%%
\begin{theorem}
Let $ (\mathcal{L},\leq , {\mathbf{G}_\mathcal{M}}) $ be an ordered complete $\:M^{\circ}\, \mathbf{G}_\mathcal{M}-M^{\bullet}S$. Suppose the mapping $\mathcal{\mathring{F}}:\mathcal{L} \longrightarrow \mathcal{L}$ with $ \eta \in [0,1)$ and $\, \gamma>0, $ satisfying the following, 
%%%%%%%%%%%%%%%%%%%%%%%%%%%%%%%%%%%%%%%%%%%%%%
\begin{equation*}
\sqrt[m]{{\mathbf{G}_\mathcal{M}}(\mathcal{\mathring{F}}\check{\nu},\mathcal{\mathring{F}}\check{\varpi},\mathcal{\mathring{F}}\check{\vartheta})} \:\: \leq \:\: \bigg[\sqrt[m] {{\mathbf{G}_\mathcal{M}}(\check{\nu},\check{\varpi},\check{\vartheta})}\: \bigg]^\mathlarger{\eta}, \tag{3.1}
\end{equation*} 	
%%%%%%%%%%%%%%%%%%%%%%%%%%%%%%%%%%%%%%%%%%%%%%
\noindent  and
%%%%%%%%%%%%%%%%%%%%%%%%%%%%%%%%%%%%%%%%%%%%%%
\begin{equation*}
{\mathbf{G}_\mathcal{M}}(\check{\nu}_{0},\mathcal{\mathring{F}}\check{\nu}_{0},\mathcal{\mathring{F}}\check{\nu}_{0}) \:\: \leq \:\: (1-\eta) \, \gamma, \tag{3.2} 	
\end{equation*}	
%%%%%%%%%%%%%%%%%%%%%%%%%%%%%%%%%%%%%%%%%%%%%%
\noindent for $ \check{\nu},\check{\varpi},\check{\vartheta} \in \overline{\mathcal{B}_{\gamma} (\check{\nu}_{0}, \gamma)}$. If for a non-increasing ($\mathbf{non-inc}$) sequence $\{\check{\nu}_{n}\} \longrightarrow s $ implies that $ s \preceq \check{\nu}_{n} $. Then, there is a point $\check{\nu}^{*}$ in $\overline{\mathcal{B}_{\gamma} (\check{\nu}_{0}, \gamma)}$ such that $ \check{\nu}^{*}=\mathcal{\mathring{F}}\check{\nu}^{*} $ and $ \mathbf{G}_\mathcal{M}(\check{\nu}^{*},\check{\nu}^{*},\check{\nu}^{*})=1. $ Furthermore, if for any two points $\check{\nu},\check{\varpi}$  in $\overline{\mathcal{B}_{\gamma} (\check{\nu}_{0}, \gamma)} $ and there exists a point $ t \in \overline{\mathcal{B}_{\gamma} (\check{\nu}_{0}, \gamma)}$ such that $\, t \preceq \check{\nu} \:$ and $\: t \preceq \check{\varpi}$, that is every two points in $ \overline{\mathcal{B}_{\gamma} (\check{\nu}_{0}, \gamma)} $ has a lower bound ($\mathbf{LB}$). Then, a point $\check{\nu}^{*} $ is unique in $ \mathcal{L}$.
	
\end{theorem}	

%%%%%%%%%%%%%%%%%%%%%%%%%%%%%%%%%%%%%%%%%%%%%%
\noindent {\textbf{Proof.}} Let $\check{\nu}_{0}$ be any arbitrary point in $\mathcal{L}$ and picard sequence $ \check{\nu}_{j+1}=\mathcal{\mathring{F}}\check{\nu}_{j} \preceq \check{\nu}_{j} \: $ for all $ n \in \mathbb{N} \cup \{0\}. $ From Ineq. (3.2), we get
%%%%%%%%%%%%%%%%%%%%%%%%%%%%%%%%%%%%%%%%%%%%%%
\begin{equation*}
{\mathbf{G}_\mathcal{M}}(\check{\nu}_{0},\check{\nu}_{1},\check{\nu}_{1}) \:\: \leq \:\: (1-\eta) \:\gamma \:\: \leq \:\: \gamma,	
\end{equation*} 
%%%%%%%%%%%%%%%%%%%%%%%%%%%%%%%%%%%%%%%%%%%%%%
\noindent implying thereby that $\: \check{\nu}_{1} \in \overline{\mathcal{B}_{\gamma} (\check{\nu}_{0}, \gamma)}. $ By multiplicative triangle inequality, we have 
%%%%%%%%%%%%%%%%%%%%%%%%%%%%%%%%%%%%%%%%%%%%%%
\begin{eqnarray*}
\sqrt[m]{{\mathbf{G}_\mathcal{M}}(\check{\nu}_{0},\check{\nu}_{2},\check{\nu}_{2})} & \leq & \sqrt[m]{{\mathbf{G}_\mathcal{M}}(\check{\nu}_{0},\check{\nu}_{1},\check{\nu}_{1})}\:.\: \sqrt[m]{{\mathbf{G}_\mathcal{M}}(\check{\nu}_{1},\check{\nu}_{2},\check{\nu}_{2})}\\[5pt]
%%%%%%%%%%%%%%%%%%%%%%%%%%%%%%%%%%%%%%%%%%%%%%
& = & \sqrt[m]{{\mathbf{G}_\mathcal{M}}(\check{\nu}_{0},\check{\nu}_{1},\check{\nu}_{1})}\:.\: \sqrt[m]{{\mathbf{G}_\mathcal{M}}(\mathcal{\mathring{F}}\check{\nu}_{0},\mathcal{\mathring{F}}\check{\nu}_{1},\mathcal{\mathring{F}}\check{\nu}_{1})}\\[5pt]
%%%%%%%%%%%%%%%%%%%%%%%%%%%%%%%%%%%%%%%%%%%%%%
& \leq & \Big[ \sqrt[m]{{\mathbf{G}_\mathcal{M}}(\check{\nu}_{0},\check{\nu}_{1},\check{\nu}_{1})}\: \Big]^\mathlarger{1+\eta}, 
\end{eqnarray*}
%%%%%%%%%%%%%%%%%%%%%%%%%%%%%%%%%%%%%%%%%%%%%%
\noindent that is,
%%%%%%%%%%%%%%%%%%%%%%%%%%%%%%%%%%%%%%%%%%%%%%
\begin{eqnarray*}
{\mathbf{G}_\mathcal{M}}(\check{\nu}_{0},\check{\nu}_{2},\check{\nu}_{2}) & \leq & \Big[ {\mathbf{G}_\mathcal{M}}(\check{\nu}_{0},\mathcal{\mathring{F}}\check{\nu}_{0},\mathcal{\mathring{F}}\check{\nu}_{0})\: \Big]^\mathlarger{1+\eta}\\[5pt]
%%%%%%%%%%%%%%%%%%%%%%%%%%%%%%%%%%%%%%%%%%%%%%
& \leq & \big[\: (1-\eta) \: \gamma\: \big]^\mathlarger{1+\eta} \:\: \leq \:\: \gamma. 
\end{eqnarray*}
%%%%%%%%%%%%%%%%%%%%%%%%%%%%%%%%%%%%%%%%%%%%%%
\noindent Then, $ \check{\nu}_{2} \in \overline{\mathcal{B}_{\gamma} (\check{\nu}_{0}, \gamma)}.\:$ Consider $\:\check{\nu}_{3}, \check{\nu}_{4}, ... , \check{\nu}_{q}\:$ for every $ q \in \mathbb{N}.\:$ Taking Ineq. (3.1) in consideration, we obtain 
%%%%%%%%%%%%%%%%%%%%%%%%%%%%%%%%%%%%%%%%%%%%%%
\begin{eqnarray*}
\sqrt[m]{{\mathbf{G}_\mathcal{M}}(\check{\nu}_{q},\check{\nu}_{q+1},\check{\nu}_{q+1})} & = & \sqrt[m]{{\mathbf{G}_\mathcal{M}}(\mathcal{\mathring{F}}\check{\nu}_{q-1},\mathcal{\mathring{F}}\check{\nu}_{q},\mathcal{\mathring{F}}\check{\nu}_{q})} \:\: \leq \:\: \Big[ \sqrt[m]{{\mathbf{G}_\mathcal{M}}(\check{\nu}_{q-1},\check{\nu}_{q},\check{\nu}_{q})}\: \Big]^\mathlarger{\eta} \notag \\[5pt]
%%%%%%%%%%%%%%%%%%%%%%%%%%%%%%%%%%%%%%%%%%%%%%
& \leq & \Big[ \sqrt[m]{{\mathbf{G}_\mathcal{M}}(\check{\nu}_{q-2},\check{\nu}_{q-1},\check{\nu}_{q-1})}\: \Big]^\mathlarger{\eta^{2}} \notag \\[5pt]
%%%%%%%%%%%%%%%%%%%%%%%%%%%%%%%%%%%%%%%%%%%%%%
& \vdots &   \notag\\
%%%%%%%%%%%%%%%%%%%%%%%%%%%%%%%%%%%%%%%%%%%%%%
& \leq & \Big[ \sqrt[m]{{\mathbf{G}_\mathcal{M}}(\check{\nu}_{0},\check{\nu}_{1},\check{\nu}_{1})}\: \Big]^\mathlarger{\eta^{\,q}}. \qquad \qquad \qquad \qquad \qquad \qquad \qquad \qquad  (3.3)
%\tag{3.3}
\end{eqnarray*}  
%%%%%%%%%%%%%%%%%%%%%%%%%%%%%%%%%%%%%%%%%%%%%%
\noindent Using Ineq. (3.1) and Ineq. (3.3), we find
%%%%%%%%%%%%%%%%%%%%%%%%%%%%%%%%%%%%%%%%%%%%%%
\begin{eqnarray*}
\sqrt[m]{{\mathbf{G}_\mathcal{M}}(\check{\nu}_{0},\check{\nu}_{q+1},\check{\nu}_{q+1})} & \leq & \sqrt[m]{{\mathbf{G}_\mathcal{M}}(\check{\nu}_{0},\check{\nu}_{1},\check{\nu}_{1})}\:.\: \sqrt[m]{{\mathbf{G}_\mathcal{M}}(\check{\nu}_{1},\check{\nu}_{2},\check{\nu}_{2})}\:.\: ... \:.\: \sqrt[m]{{\mathbf{G}_\mathcal{M}}(\check{\nu}_{q},\check{\nu}_{q+1},\check{\nu}_{q+1})}\\[5pt] 
%%%%%%%%%%%%%%%%%%%%%%%%%%%%%%%%%%%%%%%%%%%%%%
& \leq & \bigg[\sqrt[m]{{\mathbf{G}_\mathcal{M}}(\check{\nu}_{0},\check{\nu}_{1},\check{\nu}_{1})}\: \bigg]^\mathlarger{1 + \eta + ... + \eta^{\,q}},  
\end{eqnarray*}
%%%%%%%%%%%%%%%%%%%%%%%%%%%%%%%%%%%%%%%%%%%%%%
\noindent that becomes as follows 
%%%%%%%%%%%%%%%%%%%%%%%%%%%%%%%%%%%%%%%%%%%%%%
\begin{eqnarray*}
{\mathbf{G}_\mathcal{M}}(\check{\nu}_{0},\check{\nu}_{q+1},\check{\nu}_{q+1}) & \leq & \bigg[{\mathbf{G}_\mathcal{M}}(\check{\nu}_{0},\mathcal{\mathring{F}}\check{\nu}_{0},\mathcal{\mathring{F}}\check{\nu}_{0}) \bigg]^\mathlarger{\dfrac{1-\eta^{\,q+1}}{1-\eta}} \\[5pt]
%%%%%%%%%%%%%%%%%%%%%%%%%%%%%%%%%%%%%%%%%%%%%%
& \leq & \bigg[(1-\eta) \, \gamma\: \bigg]^\mathlarger{\dfrac{1-\eta^{\,q+1}}{1-\eta}} \: \leq \: \gamma.
\end{eqnarray*}
%%%%%%%%%%%%%%%%%%%%%%%%%%%%%%%%%%%%%%%%%%%%%%
\noindent Hence, $ \check{\nu}_{q+1} \in \overline{\mathcal{B}_{\gamma} (\check{\nu}_{0}, \gamma)}.\:$ Thus, $\:\check{\nu}_{j} \in \overline{\mathcal{B}_{\gamma} (\check{\nu}_{0}, \gamma)}\:$ for all $\: j\in \mathbb{N}.\:$  Consequently, Ineq. (3.3) convert to
%%%%%%%%%%%%%%%%%%%%%%%%%%%%%%%%%%%%%%%%%%%%%%
\begin{equation}
\sqrt[m]{{\mathbf{G}_\mathcal{M}}(\check{\nu}_{j},\check{\nu}_{j+1},\check{\nu}_{j+1})} \:\: \leq \:\: \Big[ \sqrt[m]{{\mathbf{G}_\mathcal{M}}(\check{\nu}_{0},\check{\nu}_{1},\check{\nu}_{1})}\: \Big]^\mathlarger{\eta^{\,j}}. \tag{3.4}
\end{equation}
%%%%%%%%%%%%%%%%%%%%%%%%%%%%%%%%%%%%%%%%%%%%%%
\noindent From Ineq. (3.4), we have  
%%%%%%%%%%%%%%%%%%%%%%%%%%%%%%%%%%%%%%%%%%%%%%
\begin{eqnarray*}
&&\sqrt[m]{{\mathbf{G}_\mathcal{M}}(\check{\nu}_{j},\check{\nu}_{j+k},\check{\nu}_{j+k})}\\[5pt] 
%%%%%%%%%%%%%%%%%%%%%%%%%%%%%%%%%%%%%%%%%%%%%%
& \leq & \sqrt[m]{{\mathbf{G}_\mathcal{M}}(\check{\nu}_{j},\check{\nu}_{j+1},\check{\nu}_{j+1})}\:.\:\sqrt[m]{{\mathbf{G}_\mathcal{M}}(\check{\nu}_{j+1},\check{\nu}_{j+2},\check{\nu}_{j+2})}\:.\: ... \:.\: \sqrt[m]{{\mathbf{G}_\mathcal{M}}(\check{\nu}_{j+k-1},\check{\nu}_{j+k},\check{\nu}_{j+k})}\\[5pt] 
%%%%%%%%%%%%%%%%%%%%%%%%%%%%%%%%%%%%%%%%%%%%%%
& \leq & \Big[ \sqrt[m]{{\mathbf{G}_\mathcal{M}}(\check{\nu}_{0},\check{\nu}_{1},\check{\nu}_{1})}\: \Big]^\mathlarger{\eta^{\,j}\:\dfrac{1-\eta^{k}}{1-\eta}} \: \longrightarrow 1, \quad \:\: j \longrightarrow +\infty.
\end{eqnarray*}
%%%%%%%%%%%%%%%%%%%%%%%%%%%%%%%%%%%%%%%%%%%%%%
\noindent This means the sequence $ \{\check{\nu}_{j}\} $ is a $\: M^{\circ}\,\mathbf{G}_\mathcal{M}-C^{\bullet} \:$ sequence in
$\: (\overline{\mathcal{B}_{\gamma} (\check{\nu}_{0},\gamma)},{\mathbf{G}_\mathcal{M}})$.\smallskip    

%%%%%%%%%%%%%%%%%%%%%%%%%%%%%%%%%%%%%%%%%%%%%%
\noindent Furthermore, there exists $\: \check{\nu}^{*} \in \overline{\mathcal{B}_{\gamma} (\check{\nu}_{0}, \gamma)} \:$ with
%%%%%%%%%%%%%%%%%%%%%%%%%%%%%%%%%%%%%%%%%%%%%%
\begin{equation}
\lim\limits_{j \longrightarrow +\infty} \sqrt[m]{\mathbf{G}_\mathcal{M}(\check{\nu}_{j},\check{\nu}^{*},\check{\nu}^{*})} \:\: = \:\: \lim\limits_{j \longrightarrow +\infty} \sqrt[m]{\mathbf{G}_\mathcal{M}(\check{\nu}^{*},\check{\nu}_{j},\check{\nu}_{j})} \:\: = \:\: 1. \tag{3.5}
\end{equation}
%%%%%%%%%%%%%%%%%%%%%%%%%%%%%%%%%%%%%%%%%%%%%%
\noindent  Now, assume that $ \check{\nu}^{*} \preceq \check{\nu}_{j} \preceq \check{\nu}_{j-1}, $ then
%%%%%%%%%%%%%%%%%%%%%%%%%%%%%%%%%%%%%%%%%%%%%%
\begin{eqnarray*}
\sqrt[m]{\mathbf{G}_\mathcal{M}(\check{\nu}^{*},\mathcal{\mathring{F}}\check{\nu}^{*},\mathcal{\mathring{F}}\check{\nu}^{*})} & \leq & \sqrt[m]{\mathbf{G}_\mathcal{M}(\check{\nu}^{*},\check{\nu}_{j},\check{\nu}_{j})} \:.\:  \sqrt[m]{\mathbf{G}_\mathcal{M}(\check{\nu}_{j},\mathcal{\mathring{F}}\check{\nu}^{*},\mathcal{\mathring{F}}\check{\nu}^{*})} \\[5pt]
%%%%%%%%%%%%%%%%%%%%%%%%%%%%%%%%%%%%%%%%%%%%%%
& = & \sqrt[m]{\mathbf{G}_\mathcal{M}(\check{\nu}^{*},\check{\nu}_{j},\check{\nu}_{j})} \:.\:  \sqrt[m]{\mathbf{G}_\mathcal{M}(\mathcal{\mathring{F}}\check{\nu}_{j-1},\mathcal{\mathring{F}}\check{\nu}^{*},\mathcal{\mathring{F}}\check{\nu}^{*})} \\[5pt]
%%%%%%%%%%%%%%%%%%%%%%%%%%%%%%%%%%%%%%%%%%%%%%
& \leq &  \sqrt[m]{\mathbf{G}_\mathcal{M}(\check{\nu}^{*},\check{\nu}_{j},\check{\nu}_{j})} \:.\:  \Big[ \sqrt[m]{\mathbf{G}_\mathcal{M}(\check{\nu}_{j-1},\check{\nu}^{*},\check{\nu}^{*})}\: \Big]^\mathlarger{\eta} \\[5pt]
%%%%%%%%%%%%%%%%%%%%%%%%%%%%%%%%%%%%%%%%%%%%%%
& \leq & \lim\limits_{j \longrightarrow +\infty} \bigg(  \sqrt[m]{\mathbf{G}_\mathcal{M}(\check{\nu}^{*},\check{\nu}_{j},\check{\nu}_{j})} \:.\:  \Big[ \sqrt[m]{\mathbf{G}_\mathcal{M}(\check{\nu}_{j-1},\check{\nu}^{*},\check{\nu}^{*})}\: \Big]^\mathlarger{\eta}\: \bigg) \:\: = \:\: 1, 
\end{eqnarray*}
%%%%%%%%%%%%%%%%%%%%%%%%%%%%%%%%%%%%%%%%%%%%%%
\noindent which is a contradiction. Then, $ \check{\nu}^{*}=\mathcal{\mathring{F}}\check{\nu}^{*}.\: $ By a similar method, $\: \mathbf{G}_\mathcal{M}(\mathcal{\mathring{F}}\check{\nu}^{*},\mathcal{\mathring{F}}\check{\nu}^{*},\check{\nu}^{*})=1\:$ and hence, $ \mathcal{\mathring{F}}\check{\nu}^{*}=\check{\nu}^{*}. $ Now,
%%%%%%%%%%%%%%%%%%%%%%%%%%%%%%%%%%%%%%%%%%%%%%
\begin{equation*}
\sqrt[m]{\mathbf{G}_\mathcal{M}(\check{\nu}^{*},\check{\nu}^{*},\check{\nu}^{*})} \: = \: \sqrt[m]{\mathbf{G}_\mathcal{M}(\mathcal{\mathring{F}}\check{\nu}^{*},\mathcal{\mathring{F}}\check{\nu}^{*},\mathcal{\mathring{F}}\check{\nu}^{*})}  \:\: \leq \:\:  \Big[\sqrt[m]{\mathbf{G}_\mathcal{M}(\check{\nu}^{*},\check{\nu}^{*},\check{\nu}^{*})}\:\Big]^\mathlarger{\eta}, 
\end{equation*}
%%%%%%%%%%%%%%%%%%%%%%%%%%%%%%%%%%%%%%%%%%%%%%
\noindent which is a contradiction, since $\eta \in [0,1).\:$ Thus, $\: \mathbf{G}_\mathcal{M}(\check{\nu}^{*},\check{\nu}^{*},\check{\nu}^{*})=1. $ \medskip

\noindent \textbf{Uniqueness:}\smallskip

\noindent Consider $ \check{\varpi}^{*} $ be another point in $ \overline{\mathcal{B}_{\gamma} (\check{\nu}_{0}, \gamma)} $ such that $ \check{\varpi}^{*} = F\check{\varpi}^{*}.\: $ If $\: \check{\nu}^{*} \:$ and $\: \check{\varpi}^{*} \:$ are comparable, then
%%%%%%%%%%%%%%%%%%%%%%%%%%%%%%%%%%%%%%%%%%%%%%
\begin{equation*}
\sqrt[m]{\mathbf{G}_\mathcal{M}(\check{\nu}^{*},\check{\varpi}^{*},\check{\varpi}^{*})} \: = \: \sqrt[m]{\mathbf{G}_\mathcal{M}(\mathcal{\mathring{F}}\check{\nu}^{*},\mathcal{\mathring{F}}\check{\varpi}^{*},\mathcal{\mathring{F}}\check{\varpi}^{*})}  \:\: \leq  \:\: \Big[\sqrt[m]{\mathbf{G}_\mathcal{M}(\check{\nu}^{*},\check{\varpi}^{*},\check{\varpi}^{*})}\:\Big]^\mathlarger{\eta}, 
\end{equation*}
%%%%%%%%%%%%%%%%%%%%%%%%%%%%%%%%%%%%%%%%%%%%%%
\noindent which is contradiction that tend us to 
$$\: \mathbf{G}_\mathcal{M}(\check{\nu}^{*},\check{\varpi}^{*},\check{\varpi}^{*})=1 \quad \text{implies}\quad \check{\nu}^{*} = \check{\varpi}^{*}.\:$$
\noindent Similarly, we can prove  $\: \mathbf{G}_\mathcal{M}(\check{\varpi}^{*},\check{\varpi}^{*},\check{\nu}^{*})=1.$ \medskip 

%%%%%%%%%%%%%%%%%%%%%%%%%%%%%%%%%%%%%%%%%%%%%%
\noindent On the other hand, If $\: \check{\nu}^{*} \:$ and $\: \check{\varpi}^{*} \:$ are not comparable then there is a point $\: t \in \overline{\mathcal{B}_{\gamma} (\check{\nu}_{0}, \gamma)} $ which is the  $\:\mathbf{LB} \:$ of $\: \check{\nu}^{*} \:$ and $\: \check{\varpi}^{*} \:$ that is $\: t \preceq \check{\nu}^{*} \:$ and $\: t \preceq \check{\varpi}^{*}. \:$ Furthermore, by argument $\: \check{\nu}^{*} \preceq \check{\nu}_{n} \:$ as $\: \check{\nu}_{n} \longrightarrow  \check{\nu}^{*}. \:$ Thus, $\: t \preceq \check{\nu}^{*} \preceq \check{\nu}_{n} \preceq ... \preceq \check{\nu}_{0}$.
%%%%%%%%%%%%%%%%%%%%%%%%%%%%%%%%%%%%%%%%%%%%%%
\begin{eqnarray*}
\sqrt[m]{\mathbf{G}_\mathcal{M}(\check{\nu}_{0},\mathcal{\mathring{F}}t,\mathcal{\mathring{F}}t)} & \leq & \sqrt[m]{\mathbf{G}_\mathcal{M}(\check{\nu}_{0},\check{\nu}_{1},\check{\nu}_{1})}\:.\: \sqrt[m]{\mathbf{G}_\mathcal{M}(\check{\nu}_{1},\mathcal{\mathring{F}}t,\mathcal{\mathring{F}}t)}\\
& = & \sqrt[m]{\mathbf{G}_\mathcal{M}(\check{\nu}_{0},\mathcal{\mathring{F}}\check{\nu}_{0},\mathcal{\mathring{F}}\check{\nu}_{0})}\:.\: \sqrt[m]{\mathbf{G}_\mathcal{M}(\mathcal{\mathring{F}}\check{\nu}_{0},\mathcal{\mathring{F}}t,\mathcal{\mathring{F}}t)}\\
& \leq & \sqrt[m]{\mathbf{G}_\mathcal{M}(\check{\nu}_{0},\mathcal{\mathring{F}}\check{\nu}_{0},\mathcal{\mathring{F}}\check{\nu}_{0})}\:.\: \Big[\sqrt[m]{\mathbf{G}_\mathcal{M}(\check{\nu}_{0},t,t)}\: \Big]^\mathlarger{\eta},  
\end{eqnarray*}
%%%%%%%%%%%%%%%%%%%%%%%%%%%%%%%%%%%%%%%%%%%%%%
\noindent that is,
%%%%%%%%%%%%%%%%%%%%%%%%%%%%%%%%%%%%%%%%%%%%%%
\begin{eqnarray*}
\mathbf{G}_\mathcal{M}(\check{\nu}_{0},\mathcal{\mathring{F}}t,\mathcal{\mathring{F}}t) & \leq & \mathbf{G}_\mathcal{M}(\check{\nu}_{0},\mathcal{\mathring{F}}\check{\nu}_{0},\mathcal{\mathring{F}}\check{\nu}_{0})\:.\: \Big[\mathbf{G}_\mathcal{M}(\check{\nu}_{0},t,t) \Big]^\mathlarger{\eta}\\
& \leq & (1-\eta)\, \gamma \:.\: \big[(1-\eta)\,\gamma \, \big]^\mathlarger{\eta}  \quad  \text{(by Ineq. (3.1) and Ineq. (3.2))}\\
& \leq & \gamma,  
\end{eqnarray*}
%%%%%%%%%%%%%%%%%%%%%%%%%%%%%%%%%%%%%%%%%%%%%%
\noindent where $\: \check{\nu}_{0}, t \in \overline{\mathcal{B}_{\gamma} (\check{\nu}_{0}, \gamma)}\: $ and this means that $\: \mathcal{\mathring{F}}t \in \overline{\mathcal{B}_{\gamma} (\check{\nu}_{0}, \gamma)}.  $ \medskip

%%%%%%%%%%%%%%%%%%%%%%%%%%%%%%%%%%%%%%%%%%%%%%
\noindent Now, we show that $ \mathcal{\mathring{F}}^{j}\,t \in \overline{\mathcal{B}_{\gamma} (\check{\nu}_{0}, \gamma)}\: $ by using mathematical induction. Suppose $ \mathcal{\mathring{F}}^{2}\,t, \mathcal{\mathring{F}}^{3}\,t, ... , \mathcal{\mathring{F}}^{q}\,t \in \overline{\mathcal{B}_{\gamma} (\check{\nu}_{0}, \gamma)}\: $ for all  $\: q \in \mathbb{N}.\: $ As $\: \mathcal{\mathring{F}}^{q}\,t \preceq \mathcal{\mathring{F}}^{q-1}\,t \preceq ... \preceq t \preceq \check{\nu}^{*} \preceq \check{\nu}_{n} \preceq ... \preceq \check{\nu}_{0},\:$ then  
%%%%%%%%%%%%%%%%%%%%%%%%%%%%%%%%%%%%%%%%%%%%%%
\begin{eqnarray*}
\sqrt[m]{\mathbf{G}_\mathcal{M}(\check{\nu}_{q+1},\mathcal{\mathring{F}}^{q+1}\,t,\mathcal{\mathring{F}}^{q+1}\,t)} & = & \sqrt[m]{\mathbf{G}_\mathcal{M}(\mathcal{\mathring{F}}\check{\nu}_{q},\mathcal{\mathring{F}}(\mathcal{\mathring{F}}^{q}\,t),\mathcal{\mathring{F}}(\mathcal{\mathring{F}}^{q}\,t))}\\
& \leq & \Big[\sqrt[m]{\mathbf{G}_\mathcal{M}(\check{\nu}_{q},\mathcal{\mathring{F}}^{q}\,t,\mathcal{\mathring{F}}^{q}\,t)}\: \Big]^\mathlarger{\eta} \leq ... \leq  \Big[\sqrt[m]{\mathbf{G}_\mathcal{M}(\check{\nu}_{q},\mathcal{\mathring{F}}^{q}\,t,\mathcal{\mathring{F}}^{q}\,t)}\: \Big]^\mathlarger{\eta^{q+1}}.  
\end{eqnarray*}  
%%%%%%%%%%%%%%%%%%%%%%%%%%%%%%%%%%%%%%%%%%%%%%
\noindent It follows that 
$$\mathbf{G}_\mathcal{M}(\check{\nu}_{q+1},\mathcal{\mathring{F}}^{q+1}\,t,\mathcal{\mathring{F}}^{q+1}\,t) \leq  \Big[\mathbf{G}_\mathcal{M}(\check{\nu}_{0},t,t) \Big]^\mathlarger{\eta^{q+1}}. \qquad (3.6)$$ 
%%%%%%%%%%%%%%%%%%%%%%%%%%%%%%%%%%%%%%%%%%%%%%
\noindent Now,  
%%%%%%%%%%%%%%%%%%%%%%%%%%%%%%%%%%%%%%%%%%%%%%
\begin{eqnarray*}
\mathbf{G}_\mathcal{M}(\check{\nu}_{0},\mathcal{\mathring{F}}^{q+1}\,t,\mathcal{\mathring{F}}^{q+1}\,t) & \leq &
\mathbf{G}_\mathcal{M}(\check{\nu}_{0},\check{\nu}_{1},\check{\nu}_{1})\:.\: ... \:.\: \mathbf{G}_\mathcal{M}(\check{\nu}_{q},\check{\nu}_{q+1},\check{\nu}_{q+1})\:.\:\mathbf{G}_\mathcal{M}(\check{\nu}_{q+1},\mathcal{\mathring{F}}^{q+1}\,t,\mathcal{\mathring{F}}^{q+1}\,t) \\
& \leq &
\mathbf{G}_\mathcal{M}(\check{\nu}_{0},\check{\nu}_{1},\check{\nu}_{1})\:.\: ... \:.\: \Big[ \mathbf{G}_\mathcal{M}(\check{\nu}_{0},\check{\nu}_{1},\check{\nu}_{1})\Big]^\mathlarger{\eta^{q}} \:.\: \Big[ \mathbf{G}_\mathcal{M}(\check{\nu}_{0},t,t)\Big]^\mathlarger{\eta^{q+1}} \\
& \leq &
\Big[ \mathbf{G}_\mathcal{M}(\check{\nu}_{0},\check{\nu}_{1},\check{\nu}_{1})\Big]^\mathlarger{1+\eta+...+\eta^{q}} \:.\: \Big[ \mathbf{G}_\mathcal{M}(\check{\nu}_{0},t,t)\Big]^\mathlarger{\eta^{q+1}} \\
& \leq &
\bigg[ (1-\eta)\, \gamma \, \bigg]^\mathlarger{\dfrac{1-\eta^{q+1}}{1-\eta}} \:.\: \Big[ (1-\eta)\,\gamma \, \Big]^\mathlarger{\eta^{q+1}}\\
& \leq & 
\bigg[ (1-\eta)\, \gamma\, \bigg]^\mathlarger{\dfrac{1-\eta^{q+2}}{1-\eta}}\:\: \leq \:\:\gamma. 
\end{eqnarray*}
%%%%%%%%%%%%%%%%%%%%%%%%%%%%%%%%%%%%%%%%%%%%%%
\noindent It means that $\: \mathcal{\mathring{F}}^{q+1}\,t \in \overline{\mathcal{B}_{\gamma} (\check{\nu}_{0}, \gamma)}\: $ and so $\: \mathcal{\mathring{F}}^{j}\,t \in \overline{\mathcal{B}_{\gamma} (\check{\nu}_{0}, \gamma)} \:$ for every $\: j \in \mathbb{N}. \:$ Further 
%%%%%%%%%%%%%%%%%%%%%%%%%%%%%%%%%%%%%%%%%%%%%%
\begin{eqnarray*}
&&\mathbf{G}_\mathcal{M}(\check{\nu}^{*},\check{\varpi}^{*},\check{\varpi}^{*}) \\
& \leq & \mathbf{G}_\mathcal{M}( \mathcal{\mathring{F}}^{j}\,\check{\nu}^{*},\mathcal{\mathring{F}}^{j-1}\,t,\mathcal{\mathring{F}}^{j-1}\,t)\:.\:\mathbf{G}_\mathcal{M}(\mathcal{\mathring{F}}^{j-1}\,t,\mathcal{\mathring{F}}^{j}\,\check{\varpi}^{*},\mathcal{\mathring{F}}^{j}\,\check{\varpi}^{*})  \\
& = & \mathbf{G}_\mathcal{M}( \mathcal{\mathring{F}}(\mathcal{\mathring{F}}^{j-1}\,\check{\nu}^{*}),\mathcal{\mathring{F}}(\mathcal{\mathring{F}}^{j-2}\,t),\mathcal{\mathring{F}}(\mathcal{\mathring{F}}^{j-2}\,t))\:.\:\mathbf{G}_\mathcal{M}(\mathcal{\mathring{F}}(\mathcal{\mathring{F}}^{j-2}\,t),\mathcal{\mathring{F}}(\mathcal{\mathring{F}}^{j-1}\,\check{\varpi}^{*}),\mathcal{\mathring{F}}(\mathcal{\mathring{F}}^{j-1}\,\check{\varpi}^{*}))\\
& \leq & \Big[ \mathbf{G}_\mathcal{M}( \mathcal{\mathring{F}}^{j-1}\check{\nu}^{*},F^{j-2} \, t, \mathcal{\mathring{F}}^{j-2} \, t) \Big]^\mathlarger{\eta}\:.\:\Big[\mathbf{G}_\mathcal{M}(\mathcal{\mathring{F}}^{j-2}\,t,\mathcal{\mathring{F}}^{j-1}\,\check{\varpi}^{*},\mathcal{\mathring{F}}^{j-1}\,\check{\varpi}^{*})\Big]^\mathlarger{\eta} \\
& \vdots & \\
&\leq & \Big[ \mathbf{G}_\mathcal{M}( \check{\nu}^{*},\mathcal{\mathring{F}}t,\mathcal{\mathring{F}}t)\Big]^\mathlarger{\eta^{j}}\:.\:\Big[\mathbf{G}_\mathcal{M}(\mathcal{\mathring{F}}t,\check{\varpi}^{*},\check{\varpi}^{*})\Big]^\mathlarger{\eta^{j}} \longrightarrow 1, \quad\:\: j\longrightarrow +\infty. 
\end{eqnarray*}
%%%%%%%%%%%%%%%%%%%%%%%%%%%%%%%%%%%%%%%%%%%%%%
\noindent Hence, $\: \mathbf{G}_\mathcal{M}(\check{\nu}^{*},\check{\varpi}^{*},\check{\varpi}^{*}) = 1 \: \Longrightarrow \: \check{\nu}^{*}=\check{\varpi}^{*}.\: $ By a similar method, we get \medskip

\qquad \qquad \qquad \qquad $\: \mathbf{G}_\mathcal{M}(\check{\varpi}^{*},\check{\varpi}^{*},\check{\nu}^{*}) = 1 \:$ implies $\: \check{\varpi}^{*}=\check{\nu}^{*}.\: $ \medskip 

\noindent Therefore, a point $\check{\nu}^{*} $ is unique in $ \mathcal{L}$.
%%%%%%%%%%%%%%%%%%%%%%%%%%%%%%%%%%%%%%%%%%%%%
\begin{corollary}
Let $(\mathcal{L},\preceq,{\mathbf{G}_\mathcal{M}})$ be an ordered complete $\:M^{\circ}\, \mathbf{G}_\mathcal{M}-M^{\bullet}S$. Suppose the mapping $\mathcal{\mathring{F}}:\mathcal{L} \longrightarrow \mathcal{L}$ with $ \eta \in [0,1)$ and $\, \gamma>0 $ satisfying the following,  	
$${\mathbf{G}_\mathcal{M}}(\mathcal{\mathring{F}}\check{\nu},\mathcal{\mathring{F}}\check{\varpi},\mathcal{\mathring{F}}\check{\vartheta})\:\: \leq \:\: \big[{\mathbf{G}_\mathcal{M}}(\check{\nu},\check{\varpi},\check{\vartheta}) \big]^\mathlarger{\eta}, \qquad (3.7)$$
\noindent  for $ \check{\nu},\check{\varpi},\check{\vartheta} \in \overline{\mathcal{B}_{\gamma} (\check{\nu}_{0}, \gamma)},$ with the condition $ (3.2) $. 

\noindent If for a $\:\mathbf{non-inc}\:$ 
sequence $\{\check{\nu}_{n}\} \longrightarrow s $ implies that $ s \preceq \check{\nu}_{n} $. Then, there is a point $\check{\nu}^{*}$ in $\overline{\mathcal{B}_{\gamma} (\check{\nu}_{0}, \gamma)}$ such that $ \check{\nu}^{*}=\mathcal{\mathring{F}}\check{\nu}^{*} $ and $ \mathbf{G}_\mathcal{M}(\check{\nu}^{*},\check{\nu}^{*},\check{\nu}^{*})=1. $ Furthermore, if for any two points $\check{\nu},\check{\varpi}$  in $\overline{\mathcal{B}_{\gamma} (\check{\nu}_{0}, \gamma)} $ then there exists a point $ t \in \overline{\mathcal{B}_{\gamma} (\check{\nu}_{0}, \gamma)}$ such that $\, t \preceq \check{\nu} \:$ and $\: t \preceq \check{\varpi}$, that is every two points in $ \overline{\mathcal{B}_{\gamma} (\check{\nu}_{0}, \gamma)} $ has a  $\:\mathbf{LB}$. Then, a point $\check{\nu}^{*} $ is unique.	
\end{corollary}
%%%%%%%%%%%%%%%%%%%%%%%%%%%%%%%%%%%%%%%%%%%%%
\begin{example}
\noindent Consider $\: \mathcal{L}=\mathbb{R^+}\cup \{0\} \:$ with $\mathbf{G}_\mathcal{M}: \mathcal{L}^{3}  \longrightarrow \mathcal{L}\:$ be a multiplicative $\: \mathbf{G}_\mathcal{M}-$metric on $\: \mathcal{L} \:$ is defined as follow:
%%%%%%%%%%%%%%%%%%%%%%%%%%%%%%%%%%%%%%%%%%%%%%
$$\mathbf{G}_\mathcal{M}(\check{\nu},\check{\varpi},\check{\vartheta}) \:\: = \:\:e^{\:\left|\check{\nu} - \check{\varpi} \right| \, + \, \left|\check{\varpi} - \check{\vartheta} \right| \, + \, \left|\check{\vartheta} - \check{\nu} \right|}.  $$
%%%%%%%%%%%%%%%%%%%%%%%%%%%%%%%%%%%%%%%%%%%%%%
\noindent Also, let $ \mathcal{\mathring{F}}:\mathcal{L} \longrightarrow \mathcal{L} $ be defined as 
\begin{equation}
\mathcal{\mathring{F}}\check{\nu} \:\: = \:\:\left\{
\begin{array}{c}	
\dfrac{\check{\nu}}{4} \qquad \:\:\:\: \text{if} \quad \check{\nu} \in \Big[\,0,\dfrac{1}{3}\,\Big);   \\[10pt]
\check{\nu} - \dfrac{1}{3} \qquad  \text{if} \quad \check{\nu} \in  \Big[\,\dfrac{1}{3},\infty\,\Big).
\end{array}
\right. \notag
\end{equation}
\noindent For $\: \check{\nu}_{0}=\dfrac{1}{3},\: \gamma=\dfrac{11}{2},\: \eta=\dfrac{5}{8}\:$ and $\: \overline{\mathcal{B}_{\gamma} (\check{\nu}_{0}, \gamma)}=\Big[\,0,\dfrac{11}{2}\,\Big],  \:$ we have 
$$(1-\eta) \, \gamma = \frac{33}{16}=  2.0625, $$
\noindent and 
%%%%%%%%%%%%%%%%%%%%%%%%%%%%%%%%%%%%%%%%%%%%%%
\begin{eqnarray*}
\mathbf{G}_\mathcal{M}(\check{\nu}_{0},\mathcal{\mathring{F}}\check{\nu}_{0},\mathcal{\mathring{F}}\check{\nu}_{0}) & = & \mathbf{G}_\mathcal{M}(\frac{1}{3},\mathcal{\mathring{F}}\frac{1}{3},\mathcal{\mathring{F}}\frac{1}{3}) \:\: = \:\: \mathbf{G}_\mathcal{M}(\frac{1}{3},0,0)\\
%%%%%%%%%%%%%%%%%%%%%%%%%%%%%%%%%%%%%%%%%%%%%%
& = & e^{2/3} =1.9477\\ 
%%%%%%%%%%%%%%%%%%%%%%%%%%%%%%%%%%%%%%%%%%%%%%
& \leq & (1-\eta)\,\gamma.  
\end{eqnarray*}

%%%%%%%%%%%%%%%%%%%%%%%%%%%%%%%%%%%%%%%%%%%%%%
\noindent \textbf{Step 1:} If $\: \check{\nu}, \check{\varpi}, \check{\vartheta} \in \Big[\,0,\dfrac{1}{3}\,\Big) \subseteq \overline{\mathcal{B}_{\gamma} (\check{\nu}_{0}, \gamma)}=\Big[\,0,\dfrac{11}{2}\,\Big],\: $ we get 
%%%%%%%%%%%%%%%%%%%%%%%%%%%%%%%%%%%%%%%%%%%%%%
\begin{eqnarray*}
\mathbf{G}_\mathcal{M}(\mathcal{\mathring{F}}\check{\nu},\mathcal{\mathring{F}}\check{\varpi},\mathcal{\mathring{F}}\check{\vartheta}) & = & e^{\:\frac{1}{4} \: ({\left|\check{\nu} - \check{\varpi} \right| \, + \, \left|\check{\varpi} - \check{\vartheta} \right| \, + \, \left|\check{\vartheta} - \check{\nu} \right|})}\\
%%%%%%%%%%%%%%%%%%%%%%%%%%%%%%%%%%%%%%%%%%%%%%
& \leq & e^{\:\frac{5}{8} \:({\left|\check{\nu} - \check{\varpi} \right| \, + \, \left|\check{\varpi} - \check{\vartheta} \right| \, + \, \left|\check{\vartheta} - \check{\nu} \right|})} \:\: = \:\: \big[{\mathbf{G}_\mathcal{M}}(x,y,z) \big]^\mathlarger{\eta}. 
\end{eqnarray*}
%%%%%%%%%%%%%%%%%%%%%%%%%%%%%%%%%%%%%%%%%%%%%%
\noindent \textbf{Step 2:} If $\: \check{\nu}, \check{\varpi}, \check{\vartheta} \in \Big[\,\dfrac{1}{3},\infty\,\Big),\: $ we have 
%%%%%%%%%%%%%%%%%%%%%%%%%%%%%%%%%%%%%%%%%%%%%%
\begin{eqnarray*}
\mathbf{G}_\mathcal{M}(\mathcal{\mathring{F}}x,\mathcal{\mathring{F}}y,\mathcal{\mathring{F}}z) & = & e^{\:\left|\check{\nu} - \check{\varpi} \right| \, + \, \left|\check{\varpi} - \check{\vartheta} \right| \, + \, \left|\check{\vartheta} - \check{\nu} \right|}\\ 
%%%%%%%%%%%%%%%%%%%%%%%%%%%%%%%%%%%%%%%%%%%%%%
& \geq & e^{\:\frac{5}{8} \: ({\left|\check{\nu} - \check{\varpi} \right| \, + \, \left|\check{\varpi} - \check{\vartheta} \right| \, + \, \left|\check{\vartheta} - \check{\nu} \right|})} \:\: = \:\: \big[{\mathbf{G}_\mathcal{M}}(\check{\nu}, \check{\varpi}, \check{\vartheta}) \big]^\mathlarger{\eta}. 
\end{eqnarray*}
%%%%%%%%%%%%%%%%%%%%%%%%%%%%%%%%%%%%%%%%%%%%%%
\noindent Clearly, the contractive condition doesn't satisfy in $ \mathcal{L} $ and is satisfied in $ \overline{\mathcal{B}_{\gamma} (\check{\nu}_{0}, \gamma)}$. Hence, all the conditions of Corollary 3.2 is verified in case of $\: \check{\nu}, \check{\varpi}, \check{\vartheta} \in \overline{\mathcal{B}_{\gamma} (\check{\nu}_{0}, \gamma)}. $  
\end{example}
%%%%%%%%%%%%%%%%%%%%%%%%%%%%%%%%%%%%%%%%%%%%%
Since every $\: M^{\circ}\, \mathbf{G}_\mathcal{M}-M^{\bullet}S \:$ generates $\:M^{\circ}MS$, we get the following corollaries.
%%%%%%%%%%%%%%%%%%%%%%%%%%%%%%%%%%%%%%%%%%%%%%
\begin{corollary}
Let $(\mathcal{L}, \preceq, d_\mathcal{M})$ be an ordered complete multiplicative $d_\mathcal{M}-$metric space ($M^{\circ}\, d_\mathcal{M}-M^{\bullet}S$). Suppose the mapping $\mathcal{\mathring{F}}:\mathcal{L} \longrightarrow \mathcal{L}$ with $ \eta \in [0,1)$ and $\, \gamma>0 $ satisfying the following,  	
$$\sqrt[m]{{d_\mathcal{M}}(\mathcal{\mathring{F}} \check{\nu},\mathcal{\mathring{F}} \check{\varpi})} \:\: \leq \:\:  \bigg[\sqrt[m]{{d_\mathcal{M}}(\check{\nu},\check{\varpi})}\: \bigg]^\mathlarger{\eta}, \qquad (3.8)$$
\noindent  and
%%%%%%%%%%%%%%%%%%%%%%%%%%%%%%%%%%%%%%%%%%%%%
$$\qquad \: {d_\mathcal{M}}(\check{\nu}_{0},\mathcal{\mathring{F}}\check{\nu}_{0}) \:\: \leq \:\:  (1-\eta) \, \gamma, \qquad \qquad \quad (3.9) $$
\noindent for $ \check{\nu}, \check{\varpi} \in \overline{\mathcal{B}_{\gamma} (\check{\nu}_{0}, \gamma)}$. If for a $\:\mathbf{non-inc}\:$ sequence $\{\check{\nu}_{n}\} \longrightarrow s $ implies that $ s \preceq \check{\nu}_{n} $. Then, there exists a point $\check{\nu}^{*}$ in $\overline{\mathcal{B}_{\gamma} (\check{\nu}_{0}, \gamma)}$ such that $ \check{\nu}^{*}=\mathcal{\mathring{F}}\check{\nu}^{*} \:$ and $\: d_\mathcal{M}(\check{\nu}^{*},\check{\nu}^{*})=1. $ Furthermore, if for any two points $\check{\nu}, \check{\varpi}$  in $\overline{\mathcal{B}_{\gamma} (\check{\nu}_{0}, \gamma)} $ then there is a point $\: t \in \overline{\mathcal{B}_{\gamma} (\check{\nu}_{0}, \gamma)}\:$ such that $\: t \preceq \check{\nu} \:$ and $\: t \preceq \check{\varpi}$, that is every two points in $\: \overline{\mathcal{B}_{\gamma} (\check{\nu}_{0}, \gamma)} \:$ has a 
$\:\mathbf{LB}$. Then, $\check{\nu}^{*} $ is a unique  point in $ \mathcal{L}$.
\end{corollary}

%%%%%%%%%%%%%%%%%%%%%%%%%%%%%%%%%%%%%%%%%%%%%
\begin{corollary}
Consider $(\mathcal{L},\preceq,d_\mathcal{M})$ be an ordered complete $\:M^{\circ}\, d_\mathcal{M}-M^{\bullet}S$. Suppose the mapping $ \mathcal{\mathring{F}}:\mathcal{L} \longrightarrow \mathcal{L}$ with $ \eta \in [0,1)$ and $\, \gamma>0 $ satisfying the following, 	
$${d_{\textit{\textbf{M}}}}(\mathcal{\mathring{F}}\check{\nu},\mathcal{\mathring{F}}\check{\varpi}) \:\: \leq \:\: \big[{d_{\textit{\textbf{M}}}}(\check{\nu}, \check{\varpi}) \big]^\mathlarger{\eta}, \qquad (3.10)$$
\noindent for $ \check{\nu}, \check{\varpi} \in \overline{\mathcal{B}_{\gamma} (\check{\nu}_{0}, \gamma)},$  with the condition $(3.9)$. 

\noindent If for a $\:\mathbf{non-inc}\:$ 
sequence $\{\check{\nu}_{n}\} \longrightarrow s \:$ implies that $\: s \preceq \check{\nu}_{n}$. Then, there is a point $\: \check{\nu}^{*}\:$ in $\: \overline{\mathcal{B}_{\gamma} (\check{\nu}_{0}, \gamma)} \:$ such that $\: \check{\nu}^{*}=\mathcal{\mathring{F}}\check{\nu}^{*} \:$ and $\: d_\mathcal{M}(\check{\nu}^{*},\check{\nu}^{*})=1$. Furthermore, if for any two points $\: \check{\nu}, \check{\varpi} \:$  in $\: \overline{\mathcal{B}_{\gamma} (\check{\nu}_{0}, \gamma)} \:$ then there exists a point $\: t \in \overline{\mathcal{B}_{\gamma} (\check{\nu}_{0}, \gamma)} \:$ such that $\: t \preceq \check{\nu} \:$ and $\: t \preceq \check{\varpi}$, that is every two points in $\: \overline{\mathcal{B}_{\gamma} (\check{\nu}_{0}, \gamma)} \:$ has a $\: \mathbf{LB}$. Then, a fixed point $ \check{\nu}^{*} $ is unique.
	
\end{corollary}

%%%%%%%%%%%%%%%%%%%%%%%%%%%%%%%%%%%%%%%%%%%%%
\begin{theorem}
Let $(\mathcal{L},\preceq,{\mathbf{G}_\mathcal{M}})$ be an ordered complete $\:M^{\circ}\, \mathbf{G}_\mathcal{M}-M^{\bullet}S$. Suppose the mapping $\mathcal{\mathring{F}}:\mathcal{L} \longrightarrow \mathcal{L}\:$  with $\: \eta \in [0,1)$ and $\, \gamma>0 $ satisfying the following, 	
$$\qquad \sqrt[m]{{\mathbf{G}_\mathcal{M}}(\mathcal{\mathring{F}}\check{\nu},\mathcal{\mathring{F}}\check{\varpi},\mathcal{\mathring{F}}\check{\vartheta})} \:\: \leq \:\: \mathcal{M}, \qquad (3.11)$$ 
\noindent since
%%%%%%%%%%%%%%%%%%%%%%%%%%%%%%%%%%%%%%%%%%%%%
\begin{equation*}
\noindent\mathcal{M} \:\: = \:\: \left[ \: \max  \, \left\{
\begin{array}{c}

\sqrt[m]{{\mathbf{G}_\mathcal{M}}(\check{\nu}, \check{\varpi}, \check{\vartheta})}, \sqrt[m]{{\mathbf{G}_\mathcal{M}}(\check{\nu}, \mathcal{\mathring{F}}\check{\nu}, \mathcal{\mathring{F}}\check{\nu})},\\[10pt]

\sqrt[m]{{\mathbf{G}_\mathcal{M}}( \check{\varpi},\mathcal{\mathring{F}} \check{\varpi}, \mathcal{\mathring{F}} \check{\varpi})},\sqrt[m]{{\mathbf{G}_\mathcal{M}}(\check{\nu},\mathcal{\mathring{F}}\check{\varpi},\mathcal{\mathring{F}}\check{\varpi})},\\[10pt]

\sqrt[m]{\min \big\{{\mathbf{G}_\mathcal{M}}(\check{\vartheta},\mathcal{\mathring{F}}\check{\nu},\mathcal{\mathring{F}}\check{\nu}),{\mathbf{G}_\mathcal{M}}(\check{\nu},\check{\vartheta},\check{\vartheta})\big\}}\\[10pt]

\end{array}
\right\}\:\right]^\mathlarger{\eta},\\[15pt] \\
\smallskip
\end{equation*} 	
\noindent  and
%%%%%%%%%%%%%%%%%%%%%%%%%%%%%%%%%%%%%%%%%%%%%
$$\qquad \qquad {\mathbf{G}_\mathcal{M}}(\check{\nu}_{0},\mathcal{\mathring{F}}\check{\nu}_{0},\mathcal{\mathring{F}}\check{\nu}_{0}) \:\: \leq \:\: (1-\eta) \, \gamma, \qquad (3.12)$$
 
\noindent for $ \: \check{\nu}, \check{\varpi}, \check{\vartheta} \in \overline{\mathcal{B}_{\gamma} (\check{\nu}_{0}, \gamma)}$.  If for a $\: \mathbf{non-inc}\:$ sequence $\: \{\check{\nu}_{n}\} \:$ in $\: \overline{\mathcal{B}_{\gamma} (\check{\nu}_{0}, \gamma)} \:$ and $\:\{\check{\nu}_{n}\} \longrightarrow v \:$ implies that $\: v \preceq \check{\nu}_{n}$. Then, there exists a unique fixed point $\: \check{\nu}^{*} \:$ such that $\: {\mathbf{G}_\mathcal{M}}(\check{\nu}^{*},\check{\nu}^{*},\check{\nu}^{*})=1\:$ and $\: \check{\nu}^{*}=\mathcal{\mathring{F}}\check{\nu}^{*}$.
	
\end{theorem}	

%%%%%%%%%%%%%%%%%%%%%%%%%%%%%%%%%%%%%%%%%%%%%
\noindent {\textbf{Proof.}} Consider an arbitrary point $\check{\nu}_0$ in $\mathcal{L}$. and a picard sequence $ \check{\nu}_{q+1}=\mathcal{\mathring{F}}\check{\nu}_{q} \preceq \check{\nu}_{q} \: $ for all $ n \in \mathbb{N} \cup \{0\}. $ From inequality (3.12), we find  
$${\mathbf{G}_\mathcal{M}}(\check{\nu}_{0},\check{\nu}_{1},\check{\nu}_{1}) \:\: \leq \:\: (1-\eta) \:\gamma \:\: \leq \:\: \gamma,$$ 
	
%%%%%%%%%%%%%%%%%%%%%%%%%%%%%%%%%%%%%%%%%%%%%
\noindent for all $ j \in \mathbb{N} \cup \{0\}$. Now, from inequalities (3.12), we obtain $ {\mathbf{G}_\mathcal{M}}(\check{\nu}_{0},\check{\nu}_{1},\check{\nu}_{1}) \leq \gamma $ and $ {\mathbf{G}_\mathcal{M}}(\check{\nu}_{1},\check{\nu}_{2},\check{\nu}_{2}) \leq \gamma, $ which tends to $ \check{\nu}_{1},\check{\nu}_{2} \in \overline{\mathcal{B}_{\gamma} (\check{\nu}_{0}, \gamma)}. $ Similarly $  \check{\nu}_{3}, ... , \check{\nu}_{q} \in \overline{\mathcal{B}_{\gamma} (\check{\nu}_{0}, \gamma)} $ for all $ q \in \mathbb{N}.\: $ Now, \medskip

\noindent $ \sqrt[m]{{\mathbf{G}_\mathcal{M}}(\check{\nu}_{q},\check{\nu}_{q+1},\check{\nu}_{q+1})} \:\: = \:\: \sqrt[m]{{\mathbf{G}_\mathcal{M}}(\mathcal{\mathring{F}}\check{\nu}_{q-1},\mathcal{\mathring{F}}\check{\nu}_{q},\mathcal{\mathring{F}}\check{\nu}_{q})} $ \medskip \medskip

$\:\:\leq \:\: \left[ \: \max \, \, \left\{
\begin{array}{c}

\sqrt[m]{{\mathbf{G}_\mathcal{M}}(\check{\nu}_{q-1},\check{\nu}_{q},\check{\nu}_{q})},\sqrt[m]{{\mathbf{G}_\mathcal{M}}(\check{\nu}_{q-1},\mathcal{\mathring{F}}\check{\nu}_{q-1},\mathcal{\mathring{F}}\check{\nu}_{q-1})},\\[10pt]

\sqrt[m]{{\mathbf{G}_\mathcal{M}}(\check{\nu}_{q},\mathcal{\mathring{F}}\check{\nu}_{q},\mathcal{\mathring{F}}\check{\nu}_{q})},\sqrt[m]{{\mathbf{G}_\mathcal{M}}(\check{\nu}_{q-1},\mathcal{\mathring{F}}\check{\nu}_{q},\mathcal{\mathring{F}}\check{\nu}_{q})}, \\[10pt]

\sqrt[m]{\min \big\{{\mathbf{G}_\mathcal{M}}(\check{\nu}_{q},\mathcal{\mathring{F}}\check{\nu}_{q-1},\mathcal{\mathring{F}}\check{\nu}_{q-1}),{\mathbf{G}_\mathcal{M}}(\check{\nu}_{q-1},\check{\nu}_{q},\check{\nu}_{q}) \big\}}\\[10pt]

\end{array}
\right\}\:\right]^\mathlarger{\eta} $ \medskip\medskip

$\:\: \leq \:\: \left[ \: \max \, \, \left\{
\begin{array}{c}

\sqrt[m]{{\mathbf{G}_\mathcal{M}}(\check{\nu}_{q-1},\check{\nu}_{q},\check{\nu}_{q})},\sqrt[m]{{\mathbf{G}_\mathcal{M}}(\check{\nu}_{q-1},\check{\nu}_{q},\check{\nu}_{q})}, \\[10pt]

\sqrt[m]{{\mathbf{G}_\mathcal{M}}(\check{\nu}_{q},\check{\nu}_{q+1},\check{\nu}_{q+1})},\sqrt[m]{{\mathbf{G}_\mathcal{M}}(\check{\nu}_{q-1},\check{\nu}_{q+1},\check{\nu}_{q+1})}, \\[10pt]

\sqrt[m]{\min \big\{{\mathbf{G}_\mathcal{M}}(\check{\nu}_{q},\check{\nu}_{q},\check{\nu}_{q}),{\mathbf{G}_\mathcal{M}}(\check{\nu}_{q-1},\check{\nu}_{q},\check{\nu}_{q})\big\}} \\[10pt]

\end{array}
\right\}\: \right] ^\mathlarger{\eta}$ \medskip\medskip 

\noindent \qquad $ \leq \:\: \left[ \: \max \, \left\{
\begin{array}{c}

\sqrt[m]{{\mathbf{G}_\mathcal{M}}(\check{\nu}_{q-1},\check{\nu}_{q},\check{\nu}_{q})},\sqrt[m]{{\mathbf{G}_\mathcal{M}}(\check{\nu}_{q-1},\check{\nu}_{q},\check{\nu}_{q})},\\[10pt]

\sqrt[m]{{\mathbf{G}_\mathcal{M}}(\check{\nu}_{q},\check{\nu}_{q+1},\check{\nu}_{q+1})}, \\[10pt]

\sqrt[m]{{\mathbf{G}_\mathcal{M}}(\check{\nu}_{q-1},\check{\nu}_{q},\check{\nu}_{q})}\:.\sqrt[m]{{\mathbf{G}_\mathcal{M}}(\check{\nu}_{q},\check{\nu}_{q+1},\check{\nu}_{q+1})}, 1 \\[10pt]

\end{array}
\right\}\:\right]^\mathlarger{\eta}, \:\: (\text{using} \,\, (G_{M_1}) \,\, \text{and} \,\, (G_{M_5}) ).$ \medskip\medskip

\noindent Implying thereby,
$$   \sqrt[m]{{\mathbf{G}_\mathcal{M}}(\check{\nu}_{q},\check{\nu}_{q+1},\check{\nu}_{q+1})} \:\: \leq \:\: \Bigg[\sqrt[m]{{\mathbf{G}_\mathcal{M}}(\check{\nu}_{q-1},\check{\nu}_{q},\check{\nu}_{q})}\:.\sqrt[m]{{\mathbf{G}_\mathcal{M}}(\check{\nu}_{q},\check{\nu}_{q+1},\check{\nu}_{q+1})}\,\Bigg]^{\mathlarger{\eta}},$$
\noindent that is 
%%%%%%%%%%%%%%%%%%%%%%%%%%%%%%%%%%%%%%%%%%%%%%
\begin{eqnarray*}
{\mathbf{G}_\mathcal{M}}(\check{\nu}_{q},\check{\nu}_{q+1},\check{\nu}_{q+1}) & \leq & \Big[{\mathbf{G}_\mathcal{M}}(\check{\nu}_{q-1},\check{\nu}_{q},\check{\nu}_{q})\Big]^\mathlarger{\mu} \\
& \leq & \Big[{\mathbf{G}_\mathcal{M}}(\check{\nu}_{q-2},\check{\nu}_{q-1},\check{\nu}_{q-1})\Big]^\mathlarger{\mu^{2}} \\
& \vdots & \\
& \leq & \Big[{\mathbf{G}_\mathcal{M}}(\check{\nu}_{0},\check{\nu}_{1},\check{\nu}_{1})\Big]^\mathlarger{\mu^{q}},
\end{eqnarray*}
%%%%%%%%%%%%%%%%%%%%%%%%%%%%%%%%%%%%%%%%%%%%%%
\noindent where $ 0 < \mu=\dfrac{\eta}{1-\eta}<\dfrac{1}{2}.\: $ Taking Ineq. (3.11) and Ineq. (3.12) in consideration, we get 
%%%%%%%%%%%%%%%%%%%%%%%%%%%%%%%%%%%%%%%%%%%%%%
\begin{eqnarray*}
{\mathbf{G}_\mathcal{M}}(\check{\nu}_{0},\check{\nu}_{q+1},\check{\nu}_{q+1}) & \leq & {\mathbf{G}_\mathcal{M}}(\check{\nu}_{0},\check{\nu}_{1},\check{\nu}_{1})\:.\:{\mathbf{G}_\mathcal{M}}(\check{\nu}_{1},\check{\nu}_{2},\check{\nu}_{2}).\: \cdots \:.\: {\mathbf{G}_\mathcal{M}}(\check{\nu}_{q},\check{\nu}_{q+1},\check{\nu}_{q+1})  \\
& \leq & \Big[{\mathbf{G}_\mathcal{M}}(\check{\nu}_{0},\check{\nu}_{1},\check{\nu}_{1})\Big]^\mathlarger{\dfrac{1-\mu^{q+1}}{1-\mu} } \\
& \leq & \Big[ (1-\eta) \, \gamma \: \Big]^\mathlarger{\dfrac{1-\mu^{q+1}}{1-\mu}} \:\: \leq  \:\: \gamma. 
\end{eqnarray*}
%%%%%%%%%%%%%%%%%%%%%%%%%%%%%%%%%%%%%%%%%%%%%%
\noindent Then, $ \check{\nu}_{q+1}  \in \overline{\mathcal{B}_{\gamma} (\check{\nu}_{0}, \gamma)}.$ Thus, $ \check{\nu}_{j} \in \overline{\mathcal{B}_{\gamma} (\check{\nu}_{0}, \gamma)} $ for every $ j \in \mathbb{N}. $ Now, Ineq. (3.13) became
%%%%%%%%%%%%%%%%%%%%%%%%%%%%%%%%%%%%%%%%%%%%%%
$${\mathbf{G}_\mathcal{M}}(\check{\nu}_{j},\check{\nu}_{j+1},\check{\nu}_{j+1})  \:\: \leq \:\: \Big[{\mathbf{G}_\mathcal{M}}(\check{\nu}_{0},\check{\nu}_{1},\check{\nu}_{1})\Big]^\mathlarger{{\mu^{j}}}. \quad (3.14)$$ 
\noindent From Ineq. (3.14), we find
%%%%%%%%%%%%%%%%%%%%%%%%%%%%%%%%%%%%%%%%%%%%%%
\begin{eqnarray*}
{\mathbf{G}_\mathcal{M}}(\check{\nu}_{j},\check{\nu}_{j+k},\check{\nu}_{j+k})   & \leq & {\mathbf{G}_\mathcal{M}}(\check{\nu}_{j},\check{\nu}_{j+1},\check{\nu}_{j+1})\:.\:{\mathbf{G}_\mathcal{M}}(\check{\nu}_{j+1},\check{\nu}_{j+2},\check{\nu}_{j+2}).\: \cdots \:.\: {\mathbf{G}_\mathcal{M}}(\check{\nu}_{j+k-1},\check{\nu}_{j+k},\check{\nu}_{j+k}) \\
& \leq & \Big[{\mathbf{G}_\mathcal{M}}(\check{\nu}_{0},\check{\nu}_{1},\check{\nu}_{1})\Big]^\mathlarger{{\mu}^{j}\:\dfrac{1-\mu^{k}}{1-\mu}}\longrightarrow 1, \quad\:\:\: j \longrightarrow +\infty. 
\end{eqnarray*}
%%%%%%%%%%%%%%%%%%%%%%%%%%%%%%%%%%%%%%%%%%%%%%
\noindent This shows that the sequence $ \{\check{\nu}_{j}\} $ is a $\: M^{\circ}\,\mathbf{G}_\mathcal{M}-C^{\bullet} \:$ sequence in
$\: (\overline{\mathcal{B}_{\gamma} (\check{\nu}_{0}, \gamma)},{\mathbf{G}_\mathcal{M}})$.  Then, there exists $\: \check{\nu}^{*} \in \overline{\mathcal{B}_{\gamma} (\check{\nu}_{0}, \gamma)} \:$ with (3.5) is verified. \medskip

\noindent Now, suppose that $\:\check{\nu}^{*} \leq \check{\nu}_{j} \leq \check{\nu}_{j-1},$ then  
%%%%%%%%%%%%%%%%%%%%%%%%%%%%%%%%%%%%%%%%%%%%%%
\begin{eqnarray*}
\sqrt[m]{\mathbf{G}_\mathcal{M}(\check{\nu}^{*},\mathcal{\mathring{F}}\check{\nu}^{*},\mathcal{\mathring{F}}\check{\nu}^{*})} & \leq & \sqrt[m]{\mathbf{G}_\mathcal{M}(\check{\nu}^{*},\check{\nu}_{j},\check{\nu}_{j})} \:.\:  \sqrt[m]{\mathbf{G}_\mathcal{M}(\check{\nu}_{j},\mathcal{\mathring{F}}\check{\nu}^{*},\mathcal{\mathring{F}}\check{\nu}^{*})} \\
& = & \sqrt[m]{\mathbf{G}_\mathcal{M}(\check{\nu}^{*},\check{\nu}_{j},\check{\nu}_{j})} \:.\:  \sqrt[m]{\mathbf{G}_\mathcal{M}(\mathcal{\mathring{F}}\check{\nu}_{j-1},\mathcal{\mathring{F}}\check{\nu}^{*},\mathcal{\mathring{F}}\check{\nu}^{*})} \\
& \leq  & \sqrt[m]{\mathbf{G}_\mathcal{M}(\check{\nu}^{*},\check{\nu}_{j},\check{\nu}_{j})} \:.\:  \Big[ \sqrt[m]{\mathbf{G}_\mathcal{M}(\check{\nu}_{j-1},\check{\nu}^{*},\check{\nu}^{*})}\: \Big]^\mathlarger{\eta} \\
& \leq & \lim\limits_{j \longrightarrow +\infty} \bigg(  \sqrt[m]{\mathbf{G}_\mathcal{M}(\check{\nu}^{*},\check{\nu}_{j},\check{\nu}_{j})} \:.\:  \Big[ \sqrt[m]{\mathbf{G}_\mathcal{M}(\check{\nu}_{j-1},\check{\nu}^{*},\check{\nu}^{*})}\: \Big]^\mathlarger{\eta}\: \bigg) \:\: = \:\: 1, 
\end{eqnarray*}
%%%%%%%%%%%%%%%%%%%%%%%%%%%%%%%%%%%%%%%%%%%%%%
\noindent which is a contradiction. Then, $ \check{\nu}^{*}=\mathcal{\mathring{F}}\check{\nu}^{*}.\: $ By a similar method, $\: \mathbf{G}_\mathcal{M}(\mathcal{\mathring{F}}\check{\nu}^{*},\mathcal{\mathring{F}}\check{\nu}^{*},\check{\nu}^{*})=1\:$ and hence $ \mathcal{\mathring{F}}\check{\nu}^{*}=\check{\nu}^{*}. $ Now,
$$ \sqrt[m]{\mathbf{G}_\mathcal{M}(\check{\nu}^{*},\check{\nu}^{*},\check{\nu}^{*})} \:\: = \:\: \sqrt[m]{\mathbf{G}_\mathcal{M}(\mathcal{\mathring{F}}\check{\nu}^{*},\mathcal{\mathring{F}}\check{\nu}^{*},\mathcal{\mathring{F}}\check{\nu}^{*})}  \:\: \leq \:\:  \Big[\sqrt[m]{\mathbf{G}_\mathcal{M}(\check{\nu}^{*},\check{\nu}^{*},\check{\nu}^{*})}\:\Big]^\mathlarger{\eta} $$ 
\noindent which is a contradiction, since $\eta \in [0,1).\:$ Thus, $\: \mathbf{G}_\mathcal{M}(\check{\nu}^{*},\check{\nu}^{*},\check{\nu}^{*})=1. $ \medskip

\noindent \textbf{Uniqueness:}\smallskip

\noindent Suppose $ \check{\varpi}^{*} $ be another point in $ \overline{\mathcal{B}_{\gamma} (\check{\nu}_{0}, \gamma)} $ such that $ \check{\varpi}^{*} = \mathcal{\mathring{F}}\check{\varpi}^{*}.\: $ If $\: \check{\nu}^{*} \:$ and $\: \check{\varpi}^{*} \:$ are comparable, then $$\sqrt[m]{\mathbf{G}_\mathcal{M}(\check{\nu}^{*},\check{\varpi}^{*},\check{\varpi}^{*})}\:\: = \:\: \sqrt[m]{\mathbf{G}_\mathcal{M}(\mathcal{\mathring{F}}\check{\nu}^{*},\mathcal{\mathring{F}}\check{\varpi}^{*},\mathcal{\mathring{F}}\check{\varpi}^{*})}  \:\: \leq  \:\: \Big[\sqrt[m]{\mathbf{G}_\mathcal{M}(\check{\nu}^{*},\check{\varpi}^{*},\check{\varpi}^{*})}\:\Big]^\mathlarger{\eta}, $$ 
\noindent which is contradiction that tend us to 
$$\: \mathbf{G}_\mathcal{M}(\check{\nu}^{*},\check{\varpi}^{*},\check{\varpi}^{*})=1 \quad \text{implies}\quad \check{\nu}^{*} = \check{\varpi}^{*}.\:$$
\noindent Similarly, we can prove  $\: \mathbf{G}_\mathcal{M}(\check{\varpi}^{*},\check{\varpi}^{*},\check{\nu}^{*})=1.$ \medskip 

\noindent On the other hand, If $\: \check{\nu}^{*} \:$ and $\: \check{\varpi}^{*} \:$ are not comparable then there exists a point $\: t \in \overline{\mathcal{B}_{\gamma} (\check{\nu}_{0}, \gamma)} $ which is the  $\:\mathbf{LB} \:$ of $\: \check{\nu}^{*} \:$ and $\: \check{\varpi}^{*} \:$ that is $\: t \preceq \check{\nu}^{*} \:$ and $\: t \preceq \check{\varpi}^{*}$. Furthermore, by argument $\: \check{\nu}^{*} \preceq \check{\nu}_{n} \:$ as $\: \check{\nu}_{n} \longrightarrow  \check{\nu}^{*}$. Thus, $\: t \preceq \check{\nu}^{*} \preceq \check{\nu}_{n} \preceq ... \preceq \check{\nu}_{0}$. Thus \medskip 
%%%%%%%%%%%%%%%%%%%%%%%%%%%%%%%%%%%%%%%%%%%%%%
\begin{eqnarray*}
\sqrt[m]{\mathbf{G}_\mathcal{M}(\check{\nu}_{0},\mathcal{\mathring{F}}t,\mathcal{\mathring{F}}t)} & \leq & \sqrt[m]{\mathbf{G}_\mathcal{M}(\check{\nu}_{0},\check{\nu}_{1},\check{\nu}_{1})}\:.\: \sqrt[m]{\mathbf{G}_\mathcal{M}(\check{\nu}_{1},\mathcal{\mathring{F}}t,\mathcal{\mathring{F}}t)} \\ 
& = & \sqrt[m]{\mathbf{G}_\mathcal{M}(\check{\nu}_{0},\mathcal{\mathring{F}}\check{\nu}_{0},\mathcal{\mathring{F}}\check{\nu}_{0})}\:.\: \sqrt[m]{\mathbf{G}_\mathcal{M}(\mathcal{\mathring{F}}\check{\nu}_{0},\mathcal{\mathring{F}}t,\mathcal{\mathring{F}}t)} \\ 
& \leq & \sqrt[m]{\mathbf{G}_\mathcal{M}(\check{\nu}_{0},\mathcal{\mathring{F}}\check{\nu}_{0},\mathcal{\mathring{F}}\check{\nu}_{0})}\:.\: \Big[\sqrt[m]{\mathbf{G}_\mathcal{M}(\check{\nu}_{0},t,t)}\: \Big]^\mathlarger{\eta},  
\end{eqnarray*}
%%%%%%%%%%%%%%%%%%%%%%%%%%%%%%%%%%%%%%%%%%%%%%
\noindent that is
%%%%%%%%%%%%%%%%%%%%%%%%%%%%%%%%%%%%%%%%%%%%%%
\begin{eqnarray*}
\mathbf{G}_\mathcal{M}(\check{\nu}_{0},\mathcal{\mathring{F}}t,\mathcal{\mathring{F}}t) & \leq & \mathbf{G}_\mathcal{M}(\check{\nu}_{0},\mathcal{\mathring{F}}\check{\nu}_{0},\mathcal{\mathring{F}}\check{\nu}_{0})\:.\: \Big[\mathbf{G}_\mathcal{M}(\check{\nu}_{0},t,t) \Big]^\mathlarger{\eta} \\
& \leq & (1-\eta)\, \gamma \:.\: \big[(1-\eta)\,\gamma \, \big]^\mathlarger{\eta}  \qquad  \text{(by Ineq. (3.12))} \\
& \leq & \gamma,  
\end{eqnarray*}
%%%%%%%%%%%%%%%%%%%%%%%%%%%%%%%%%%%%%%%%%%%%%%
\noindent where $\: \check{\nu}_{0}, t \in \overline{\mathcal{B}_{\gamma} (\check{\nu}_{0}, \gamma)}\: $ and this means that $\: \mathcal{\mathring{F}}t \in \overline{\mathcal{B}_{\gamma} (\check{\nu}_{0}, \gamma)}.  $ \medskip

\noindent Now, we prove that $ \mathcal{\mathring{F}}^{j}t \in \overline{\mathcal{B}_{\gamma} (\check{\nu}_{0}, \gamma)}\: $ by using mathematical induction. Suppose $ \mathcal{\mathring{F}}^{2}t, \mathcal{\mathring{F}}^{3}t, ... , \mathcal{\mathring{F}}^{q}t \in \overline{\mathcal{B}_{\gamma} (\check{\nu}_{0}, \gamma)}\: $ for all  $\: q \in \mathbb{N}.\: $ As $\: \mathcal{\mathring{F}}^{q}t \preceq \mathcal{\mathring{F}}^{q-1}t \preceq ... \preceq t \preceq \check{\nu}^{*} \preceq \check{\nu}_{n} \preceq ... \preceq \check{\nu}_{0},\:$ then,  
%%%%%%%%%%%%%%%%%%%%%%%%%%%%%%%%%%%%%%%%%%%%%%
\begin{eqnarray*}
\sqrt[m]{\mathbf{G}_\mathcal{M}(\check{\nu}_{q+1},\mathcal{\mathring{F}}^{q+1}t,\mathcal{\mathring{F}}^{q+1}t)} & = & \sqrt[m]{\mathbf{G}_\mathcal{M}(\mathcal{\mathring{F}}\check{\nu}_{q},\mathcal{\mathring{F}}(\mathcal{\mathring{F}}^{q}t),\mathcal{\mathring{F}}(\mathcal{\mathring{F}}^{q}t))} \\
& \leq & \Big[\sqrt[m]{\mathbf{G}_\mathcal{M}(\check{\nu}_{q},\mathcal{\mathring{F}}^{q}t,\mathcal{\mathring{F}}^{q}t)}\: \Big]^\mathlarger{\eta} \:\: \leq \:\: ... \:\: \leq \:\:  \Big[\sqrt[m]{\mathbf{G}_\mathcal{M}(\check{\nu}_{q},\mathcal{\mathring{F}}^{q}t,\mathcal{\mathring{F}}^{q}t)}\: \Big]^\mathlarger{\eta^{q+1}}. 
\end{eqnarray*}
%%%%%%%%%%%%%%%%%%%%%%%%%%%%%%%%%%%%%%%%%%%%%%
\noindent It follows that 
$$\mathbf{G}_\mathcal{M}(\check{\nu}_{q+1},\mathcal{\mathring{F}}^{q+1}t,\mathcal{\mathring{F}}^{q+1}t) \:\: \leq  \:\: \Big[\mathbf{G}_\mathcal{M}(\check{\nu}_{0},t,t) \Big]^\mathlarger{\eta^{q+1}}. \qquad (3.15)$$ 
\noindent Now, 
%%%%%%%%%%%%%%%%%%%%%%%%%%%%%%%%%%%%%%%%%%%%%%
\begin{eqnarray*}
\mathbf{G}_\mathcal{M}(\check{\nu}_{0},\mathcal{\mathring{F}}^{q+1}t,\mathcal{\mathring{F}}^{q+1}t) & \leq &
\mathbf{G}_\mathcal{M}(\check{\nu}_{0},\check{\nu}_{1},\check{\nu}_{1})\:.\: ... \:.\: \mathbf{G}_\mathcal{M}(\check{\nu}_{q},\check{\nu}_{q+1},\check{\nu}_{q+1})\:.\:\mathbf{G}_\mathcal{M}(\check{\nu}_{q+1},\mathcal{\mathring{F}}^{q+1}t,\mathcal{\mathring{F}}^{q+1}t) \\
& \leq &
\mathbf{G}_\mathcal{M}(\check{\nu}_{0},\check{\nu}_{1},\check{\nu}_{1})\:.\: ... \:.\: \Big[ \mathbf{G}_\mathcal{M}(\check{\nu}_{0},\check{\nu}_{1},\check{\nu}_{1})\Big]^\mathlarger{\eta^{q}} \:.\: \Big[ \mathbf{G}_\mathcal{M}(\check{\nu}_{0},t,t)\Big]^\mathlarger{\eta^{q+1}} \\
& \leq &
\Big[ \mathbf{G}_\mathcal{M}(\check{\nu}_{0},\check{\nu}_{1},\check{\nu}_{1})\Big]^\mathlarger{1+\eta+...+\eta^{q}} \:.\: \Big[ \mathbf{G}_\mathcal{M}(\check{\nu}_{0},t,t)\Big]^\mathlarger{\eta^{q+1}} \\
& \leq &
\bigg[ (1-\eta)\, \gamma \, \bigg]^\mathlarger{\dfrac{1-\eta^{q+1}}{1-\eta}} \:.\: \Big[ (1-\eta)\,\gamma \, \Big]^\mathlarger{\eta^{q+1}} \\
&  \leq &
\bigg[ (1-\eta)\, \gamma\, \bigg]^\mathlarger{\dfrac{1-\eta^{q+2}}{1-\eta}} \:\: \leq \:\: \gamma. 
\end{eqnarray*}
%%%%%%%%%%%%%%%%%%%%%%%%%%%%%%%%%%%%%%%%%%%%%%
\noindent It follows that $\: \mathcal{\mathring{F}}^{q+1}t \in \overline{\mathcal{B}_{\gamma} (\check{\nu}_{0}, \gamma)}\: $ and so $\: \mathcal{\mathring{F}}^{j}t \in \overline{\mathcal{B}_{\gamma} (\check{\nu}_{0}, \gamma)} \:$ for every $\: j \in \mathbb{N}. \:$ Furthermore 
%%%%%%%%%%%%%%%%%%%%%%%%%%%%%%%%%%%%%%%%%%%%%%
\begin{eqnarray*}
\mathbf{G}_\mathcal{M}(\check{\nu}^{*},\check{\varpi}^{*},\check{\varpi}^{*}) & \leq & \mathbf{G}_\mathcal{M}( \mathcal{\mathring{F}}^{j}\check{\nu}^{*},\mathcal{\mathring{F}}^{j-1}t,\mathcal{\mathring{F}}^{j-1}t)\:.\:\mathbf{G}_\mathcal{M}(\mathcal{\mathring{F}}^{j-1}t,\mathcal{\mathring{F}}^{j}\check{\varpi}^{*},\mathcal{\mathring{F}}^{j}\check{\varpi}^{*})  \\
& = & \mathbf{G}_\mathcal{M}( \mathcal{\mathring{F}}(\mathcal{\mathring{F}}^{j-1}\check{\nu}^{*}),\mathcal{\mathring{F}}(\mathcal{\mathring{F}}^{j-2}t),\mathcal{\mathring{F}}(\mathcal{\mathring{F}}^{j-2}t))\:.\:\mathbf{G}_\mathcal{M}(\mathcal{\mathring{F}}(\mathcal{\mathring{F}}^{j-2}t),\mathcal{\mathring{F}}(\mathcal{\mathring{F}}^{j-1}\check{\varpi}^{*}),\mathcal{\mathring{F}}(\mathcal{\mathring{F}}^{j-1}\check{\varpi}^{*})) \\
& \leq & \Big[ \mathbf{G}_\mathcal{M}( \mathcal{\mathring{F}}^{j-1}\check{\nu}^{*},\mathcal{\mathring{F}}^{j-2}t,\mathcal{\mathring{F}}^{j-2}t)\Big]^\mathlarger{\eta}\:.\:\Big[\mathbf{G}_\mathcal{M}(\mathcal{\mathring{F}}^{j-2}t,\mathcal{\mathring{F}}^{j-1}\check{\varpi}^{*},\mathcal{\mathring{F}}^{j-1}\check{\varpi}^{*})\Big]^\mathlarger{\eta} \\
& \vdots & \\
& \leq & \Big[ \mathbf{G}_\mathcal{M}( \check{\nu}^{*},\mathcal{\mathring{F}}t,\mathcal{\mathring{F}}t)\Big]^\mathlarger{\eta^{j}}\:.\:\Big[\mathbf{G}_\mathcal{M}(\mathcal{\mathring{F}}t,\check{\varpi}^{*},\check{\varpi}^{*})\Big]^\mathlarger{\eta^{j}} \longrightarrow 1, \quad\:\: j\longrightarrow +\infty.  
\end{eqnarray*}
%%%%%%%%%%%%%%%%%%%%%%%%%%%%%%%%%%%%%%%%%%%%%%
\noindent Hence, $\: \mathbf{G}_\mathcal{M}(\check{\nu}^{*},\check{\varpi}^{*},\check{\varpi}^{*}) = 1 \: \Longrightarrow \: \check{\nu}^{*}=\check{\varpi}^{*}.\: $ Similarly, \medskip

\qquad \qquad \qquad \qquad $\: \mathbf{G}_\mathcal{M}(\check{\varpi}^{*},\check{\varpi}^{*},\check{\nu}^{*}) = 1 \:$ implies $\: \check{\varpi}^{*}=\check{\nu}^{*}.\: $ \medskip 

\noindent Therefore, a point $\check{\nu}^{*} $ is unique in $ \mathcal{L}$. \medskip 

As illustrated, Theorem 3.1  considers a corollary to Theorem 3.6.

%%%%%%%%%%%%%%%%%%%%%%%%%%%%%%%%%%%%%%%%%%%%%
\begin{example}
\noindent Consider $\: \mathcal{L}=\mathbb{R^+}\cup \{0\} \:$ with $\mathbf{G}_\mathcal{M}: \mathcal{L}^{3}  \longrightarrow \mathcal{L}\:$ be a multiplicative $\: \mathbf{G}_\mathcal{M}-$metric on $\: \mathcal{L} \:$ is defined by
$$\mathbf{G}_\mathcal{M}(\check{\nu},\check{\varpi},\check{\vartheta}) \:\: = \:\: e^{\:\left|\check{\nu} - \check{\varpi} \right| \, + \, \left|\check{\varpi} - \check{\vartheta} \right| \, + \, \left|\check{\vartheta} - \check{\nu} \right|}.  $$
\noindent Also, let the mapping $ \mathcal{\mathring{F}}:\mathcal{L} \longrightarrow \mathcal{L} $ be defined as 
\begin{equation}
\mathcal{\mathring{F}}\check{\nu} \:\: = \:\: \left\{
\begin{array}{c}	
\dfrac{\check{\nu}}{2} \qquad \:\:\:\: \text{if} \quad \check{\nu} \in \Big(\,0,\dfrac{1}{2}\,\Big)  \cap \mathcal{L}; \\[10pt]
\check{\nu} - \dfrac{1}{4} \qquad  \text{if} \quad \check{\nu} \in  \Big[\,\dfrac{1}{2},\infty\,\Big) \cap \mathcal{L},
\end{array}
\right. \notag
\end{equation}
\noindent and
%%%%%%%%%%%%%%%%%%%%%%%%%%%%%%%%%%%%%%%%%%%%%%
\begin{eqnarray*}
\mathcal{M} & = & \left[ \: \max  \, \left\{
\begin{array}{c}
	
\sqrt[m]{{\mathbf{G}_\mathcal{M}}(\check{\nu}, \check{\varpi}, \check{\vartheta})}, \sqrt[m]{{\mathbf{G}_\mathcal{M}}(\check{\nu}, \mathcal{\mathring{F}}\check{\nu}, \mathcal{\mathring{F}}\check{\nu})},\\[10pt]
	
\sqrt[m]{{\mathbf{G}_\mathcal{M}}( \check{\varpi},\mathcal{\mathring{F}} \check{\varpi}, \mathcal{\mathring{F}} \check{\varpi})},\sqrt[m]{{\mathbf{G}_\mathcal{M}}(\check{\nu},\mathcal{\mathring{F}}\check{\varpi},\mathcal{\mathring{F}}\check{\varpi})},\\[10pt]
	
\sqrt[m]{\min \big\{{\mathbf{G}_\mathcal{M}}(\check{\vartheta},\mathcal{\mathring{F}}\check{\nu},\mathcal{\mathring{F}}\check{\nu}),{\mathbf{G}_\mathcal{M}}(\check{\nu},\check{\vartheta},\check{\vartheta})\big\}}\\[10pt]	
\end{array}
\right\}\:\right]^\mathlarger{\eta}.
\end{eqnarray*}\smallskip

%%%%%%%%%%%%%%%%%%%%%%%%%%%%%%%%%%%%%%%%%%%%%%
\noindent For $ \check{\nu}_{0}=\dfrac{1}{3},\: \gamma=\dfrac{11}{2},\: \eta=\dfrac{5}{8}\:$ and $\: \overline{\mathcal{B}_{\gamma} (\check{\nu}_{0}, \gamma)}=\Big[\,0,\dfrac{11}{2}\,\Big],  \:$ we have 
$$(1-\eta) \, \gamma \:\: = \:\: \dfrac{33}{16} \:\: = \:\:  2.0625, $$
\noindent and 
%%%%%%%%%%%%%%%%%%%%%%%%%%%%%%%%%%%%%%%%%%%%%%
\begin{eqnarray*}
\mathbf{G}_\mathcal{M}(\check{\nu}_{0},\mathcal{\mathring{F}}\check{\nu}_{0},\mathcal{\mathring{F}}\check{\nu}_{0}) & = & \mathbf{G}_\mathcal{M}(\dfrac{1}{3},\mathcal{\mathring{F}}\frac{1}{3},\mathcal{\mathring{F}}\frac{1}{3}) \:\: = \:\: \mathbf{G}_\mathcal{M}(\dfrac{1}{3},\dfrac{1}{6},\dfrac{1}{6})\\
%%%%%%%%%%%%%%%%%%%%%%%%%%%%%%%%%%%%%%%%%%%%%%
& = & e^{1/3} \:\: = \:\: 1.3956 \\ 
%%%%%%%%%%%%%%%%%%%%%%%%%%%%%%%%%%%%%%%%%%%%%%
& \leq & (1-\eta)\,\gamma.  
\end{eqnarray*}
%%%%%%%%%%%%%%%%%%%%%%%%%%%%%%%%%%%%%%%%%%%%%%
\noindent \textbf{Step 1:} If $\: \check{\nu}, \check{\varpi}, \check{\vartheta} \in \Big(\:0,\dfrac{1}{2}\:\Big) \cap \mathcal{L} \subseteq \overline{\mathcal{B}_{\gamma} (\check{\nu}_{0}, \gamma)}=\Big[\:0,\dfrac{11}{2}\:\Big],\: $ we obtain
%%%%%%%%%%%%%%%%%%%%%%%%%%%%%%%%%%%%%%%%%%%%%%
\begin{eqnarray*}
\sqrt[m]{ \mathbf{G}_\mathcal{M}(\mathcal{\mathring{F}}\check{\nu},\mathcal{\mathring{F}}\check{\varpi},\mathcal{\mathring{F}}\check{\vartheta})} & = & \sqrt[m]{ e^{\:\frac{1}{2} \: ({\left|\check{\nu} - \check{\varpi} \right| \, + \, \left|\check{\varpi} - \check{\vartheta} \right| \, + \, \left|\check{\vartheta} - \check{\nu} \right|})}}\\[7pt]
%%%%%%%%%%%%%%%%%%%%%%%%%%%%%%%%%%%%%%%%%%%%%%
& \leq & \left[ \: \max  \, \left\{
\begin{array}{c}
	
\sqrt[m]{e^{\:\left|\check{\nu} - \check{\varpi} \right| \, + \, \left|\check{\varpi} - \check{\vartheta} \right| \, + \, \left|\check{\vartheta} - \check{\nu} \right|\:}}, \sqrt[m]{e^{\:\left| \check{\nu} \right|\:}},\\[10pt]
	
\sqrt[m]{e^{\:\left| \check{\varpi} \right|\:}},\sqrt[m]{e^{\:\left| \check{\varpi} - 2 \check{\nu} \right|\:}},\\[10pt]
	
\sqrt[m]{\min \big\{e^{\:\left| \check{\nu} - 2 \check{\vartheta} \right|\:}, e^{\:2 \, \left| \check{\nu} - \check{\vartheta} \right|\:}\big\}}\\[10pt]	
\end{array}
\right\}\:\right]^\mathlarger{\eta}\\[7pt]
%%%%%%%%%%%%%%%%%%%%%%%%%%%%%%%%%%%%%%%%%%%%%%
& = & \Bigg[\sqrt[m]{e^{\:\left|\check{\nu} - \check{\varpi} \right| \, + \, \left|\check{\varpi} - \check{\vartheta} \right| \, + \, \left|\check{\vartheta} - \check{\nu} \right|\:}} \:\: \Bigg]^{\mathlarger{\eta}} \\[7pt]
& = & \Bigg[\sqrt[m]{\mathbf{G}_\mathcal{M}(\check{\nu},\check{\varpi},\check{\vartheta})} \,\, \Bigg]^{\mathlarger{\eta}}. 
\end{eqnarray*}
%%%%%%%%%%%%%%%%%%%%%%%%%%%%%%%%%%%%%%%%%%%%%%
\noindent\textbf{Step 2:} If $\: \check{\nu}, \check{\varpi}, \check{\vartheta} \in \Big[\,\dfrac{1}{2},\infty\,\Big) \cap \mathcal{L},\: $ we have 
%%%%%%%%%%%%%%%%%%%%%%%%%%%%%%%%%%%%%%%%%%%%%%
\begin{eqnarray*}
\sqrt[m]{ \mathbf{G}_\mathcal{M}(\mathcal{\mathring{F}}x,\mathcal{\mathring{F}}y,\mathcal{\mathring{F}}z) \,} & = & \sqrt[m]{ e^{\:\left|\check{\nu} - \check{\varpi} \right| \, + \, \left|\check{\varpi} - \check{\vartheta} \right| \, + \, \left|\check{\vartheta} - \check{\nu} \right|}\,}\\ 
%%%%%%%%%%%%%%%%%%%%%%%%%%%%%%%%%%%%%%%%%%%%%%
& \geq & \left[ \: \max  \, \left\{
\begin{array}{c}
	
\sqrt[m]{e^{\:\left|\check{\nu} - \check{\varpi} \right| \, + \, \left|\check{\varpi} - \check{\vartheta} \right| \, + \, \left|\check{\vartheta} - \check{\nu} \right|\:}}, \sqrt[m]{e^{\frac{ 1}{2}\:}},\\[10pt]
	
\sqrt[m]{e^{\frac{1}{2}\:}},\sqrt[m]{e^{\:2 \: \left| \check{\nu} \, - \, \check{\varpi} \, + \, \frac{1}{4} \right|\:}},\\[10pt]
	
\sqrt[m]{\min \big\{e^{\: 2 \:\left| \check{\vartheta} \, - \, \check{\nu}  \, +\, \frac{1}{4} \right|\:}, e^{\:2 \, \left| \check{\nu} - \check{\vartheta} \right|\:}\big\}}\\[10pt]	
\end{array}
\right\}\:\right]^\mathlarger{\eta}\\[7pt]
%%%%%%%%%%%%%%%%%%%%%%%%%%%%%%%%%%%%%%%%%%%%%%
& = & \Bigg[\sqrt[m]{e^{\:\left|\check{\nu} - \check{\varpi} \right| \, + \, \left|\check{\varpi} - \check{\vartheta} \right| \, + \, \left|\check{\vartheta} - \check{\nu} \right|\:}} \:\: \Bigg]^{\mathlarger{\eta}} \\[7pt]
%%%%%%%%%%%%%%%%%%%%%%%%%%%%%%%%%%%%%%%%%%%%%%
& = & \Bigg[\sqrt[m]{e^{\:\left|\check{\nu} - \check{\varpi} \right| \, + \, \left|\check{\varpi} - \check{\vartheta} \right| \, + \, \left|\check{\vartheta} - \check{\nu} \right|\:}} \:\: \Bigg]^{\mathlarger{5/8}}.
\end{eqnarray*}
%%%%%%%%%%%%%%%%%%%%%%%%%%%%%%%%%%%%%%%%%%%%%%
\noindent Clearly, the contractive condition doesn't verify in $ \Big[\,\dfrac{1}{2},\infty\,\Big) \cap\mathcal{L} $ and is verified in $ \overline{\mathcal{B}_{\gamma} (\check{\nu}_{0}, \gamma)}$. Hence, all the assertions of Theorem 3.6 is satisfied in case of $\: \check{\nu}, \check{\varpi}, \check{\vartheta} \in \overline{\mathcal{B}_{\gamma} (\check{\nu}_{0}, \gamma)}$.  
\end{example}
\newpage
\section{Conclusions}
\indent \indent In this article, we achieve some fixed point results satisfying a generalized $\Delta$-implicit contractive conditions in the context of ordered complete multiplicative $ \mathbf{G}_\mathcal{M}-$metric space ($\:M^{\circ}\, \mathbf{G}_\mathcal{M}-M^{\bullet}S$). Our results are considered a generalization and extension to the results in the literature. Some new definitons and examples are introduced in such spaces. Moreover, some examples are given to validate our obtained new results. 

%%%%%%%%%%%%%%%%%%%%%%%%%%%%%%%%%%%%%%%%%%%%%
%\section*{Conflict of interest}
%\indent \indent The authors declare that they have no conflict of interest.

%%%%%%%%%%%%%%%%%%%%%%%%%%%%%%%%%%%%%%%%%%%%%
%\section*{Acknowledgements}
%
%The authors would like to express of their great thanks to editor-in-chief and the referees for their helpful and valuable comments and suggestions.

%%%%%%%%%%%%%%%%%%%%%%%%%%%%%%%%%%%%%%%%%%%%%%%%%%%%%%%%%%%%%%%%%%%%%%%%%%%%%%%%%%%%%%%%%%%%

\end{document}